\documentclass[pdflatex,sn-mathphys-num]{sn-jnl}

\usepackage{graphicx}%
\usepackage{multirow}%
\usepackage{amsmath,amssymb,amsfonts}%
\usepackage{amsthm}%
\usepackage[title]{appendix}%
\usepackage{xcolor}%
\usepackage{textcomp}%
\usepackage{manyfoot}%
\usepackage{booktabs}%
\usepackage{algorithm}%
\usepackage{algorithmicx}%
\usepackage{algpseudocode}%
\usepackage{listings}%
\usepackage{mathtools}%
\usepackage{placeins}%

\theoremstyle{plain}
\newtheorem{theorem}{Theorem}[section]
\newtheorem{proposition}{Proposition}[section]%
\newtheorem{lemma}{Lemma}[section]%
\newtheorem{corollary}{Corollary}[section]%

\theoremstyle{remark}
\newtheorem{remark}{Remark}[section]%

\theoremstyle{definition}
\newtheorem{definition}{Definition}[section]%
\newtheorem{assumption}{Assumption}%

\counterwithin{algorithm}{section}%
\makeatletter
\AtBeginDocument{%
  \renewenvironment{proof}[1][\proofname]{\par\removelastskip\pushQED{\qed}%
    \normalfont\topsep7.5\p@\@plus7.5\p@\relax\trivlist
    \item[\hskip\labelsep\itshape#1\@addpunct{.}]\ignorespaces}%
  {\popQED\endtrivlist\@endpefalse}}
\makeatother

\begin{document}

\title[Block coordinate descent for multiobjective optimization]{A block coordinate descent method for nonconvex multiobjective optimization under H\"older continuity}





\author*[1]{\fnm{Lizet} \sur{Santa Cruz Calderón}}\email{lizet.calderon@ufms.br}

\author[2]{\fnm{John Frank} \sur{Matos Ascona}}\email{john.ascona@uems.br}
\equalcont{These authors contributed equally to this work.}

\author[3]{\fnm{Vitaliano de Sousa} \sur{Amaral}}\email{vitalianoamaral@ufpi.edu.br}
\equalcont{These authors contributed equally to this work.}

\author[4]{\fnm{Danilo Rodrigues} \spfx{de} \sur{Souza}}\email{souzadr@ime.unicamp.br}
\equalcont{These authors contributed equally to this work.}




\affil*[1]{\orgname{Federal University of Mato Grosso do Sul}, \orgaddress{\street{Rua Oscar Trindade de Barros, 740}, \city{Aquidauana}, \postcode{79200-000}, \state{MS}, \country{Brazil}}}

\affil[2]{\orgname{State University of Mato Grosso do Sul}, \orgaddress{\street{Rodovia Dourados-Itahum, km 12}, \city{Dourados}, \postcode{79804-970}, \state{MS}, \country{Brazil}}}

\affil[3]{\orgdiv{Department of Mathematics}, \orgname{Federal University of Piauí}, \orgaddress{\street{Avenida Petrônio Portela}, \city{Teresina}, \postcode{64049-550}, \state{PI}, \country{Brazil}}}

\affil[4]{\orgdiv{Department of Applied Mathematics}, \orgname{Institute of Mathematics, Statistics and Scientific Computing, University of Campinas}, \orgaddress{\street{Rua Sérgio Buarque de Holanda}, \city{Campinas}, \postcode{13083-859}, \state{SP}, \country{Brazil}}}



\abstract{In this paper, we propose a block coordinate descent method for unconstrained multiobjective optimization problems with possibly nonconvex objectives and H\"older continuous gradients, intended for large-scale instances, in which updating all the variables simultaneously is expensive. In each iteration, only one block is updated along a descent direction common to all objectives, obtained by minimizing the maximum of their regularized quadratic models. We show that, unlike the single-objective case, stationarity on each block does not imply Pareto stationarity, and that no method that updates only the blocks of a partition of the variables and is monotone in all the objectives reaches, in general, an $\epsilon$-Pareto stationary point. Therefore, when none of the blockwise stationarity measures exceeds the tolerance, the algorithm performs a joint step across all coordinates and terminates only at $\epsilon$-Pareto stationary points. We establish a complexity bound of $O(\epsilon^{-(\beta+1)/\beta})$ iterations, where $\beta$ is the smallest H\"older exponent, the same order as that of methods that update all the variables, as well as bounds on the numbers of function and gradient evaluations, which also hold with one regularization parameter per block.
}

\keywords{ Multiobjective optimization; block coordinate descent; nonconvex optimization; Hölder continuity; worst-case complexity; Pareto stationarity.}

\pacs[MSC Classification]{90C29, 90C26, 90C30, 90C06, 90C60, 65K05.}

\maketitle
\section{Introduction}\label{sec1}

This paper deals with unconstrained multiobjective problems with many variables and with nonconvex objectives whose gradients are only H\"older continuous. A point is Pareto optimal when none of the objectives can be decreased without increasing another one, and Pareto stationarity is the corresponding first-order necessary condition. Problems of this kind arise, for instance, when several learning tasks share the same parameters \cite{SK18}. Replacing the objectives by a weighted sum turns the problem into a scalar one, but the weights must be chosen beforehand and, without convexity, some Pareto optimal points cannot be obtained in this way \cite{das1997}. We work instead with descent methods that handle all the objectives together and do not require weights.

In the steepest descent method of Fliege and Svaiter \cite{FS00}, the search direction solves a min-max subproblem built from the gradients of all the objectives, and every accumulation point of the iterates is Pareto stationary. The same subproblem underlies the steepest descent method for vector optimization \cite{GS05} and the multiple-gradient descent algorithm \cite{Des12}, whose basic version, without rescaling of the gradients, yields the same direction and was used to train neural networks with several tasks \cite{SK18}. Newton \cite{FGS09}, quasi-Newton \cite{Pov14}, nonlinear conjugate gradient \cite{LP18}, Barzilai--Borwein \cite{morovati2016,chen2023} and proximal gradient \cite{TFY19} methods keep this structure and change the model or the step size.

The worst-case complexity of these methods follows the scalar theory. For nonconvex objectives with Lipschitz continuous gradients, the multiobjective steepest descent method reaches an $\epsilon$-Pareto stationary point in $O(\epsilon^{-2})$ iterations \cite{FVV19}, and Sampaio \cite{Sam26} studied lower complexity bounds for first-order methods. Calderón, Diniz-Ehrhardt and Martínez \cite{CDM20} obtained complexity bounds for high-order regularization methods under H\"older conditions on the derivatives of order $p$. H\"older conditions were also used in a partially derivative-free proximal method for composite problems \cite{AAS25}, in generalized conditional gradient methods \cite{CTY25} and in a universal nonmonotone line search method that does not need the H\"older constants \cite{PG25}. All these methods update every variable at every iteration, so the cost of an iteration grows with the full dimension of the problem.

In scalar optimization, a standard way of reducing this cost is to update one block of variables at a time \cite{Tse01,BT13,Wri15}. For convex problems with random block selection, the complexity of such methods is comparable to that of full gradient methods, while each iteration is much cheaper \cite{Nes12}. Without convexity, the price of working on blocks appears in the bounds of Amaral, Andreani, Birgin, Marcondes and Martínez \cite{AABMM22} for high-order coordinate descent with box constraints: approximate stationarity with respect to the variables of each block is reached in $O(\epsilon^{-(p+1)/p})$ iterations, whereas stationarity with respect to all the variables requires $O(\epsilon^{-(p+1)})$ iterations. Birgin and Martínez \cite{BM22} analyzed block methods with quadratic regularization for constrained problems under Lipschitz conditions.

H\"older continuity of the gradients is weaker than Lipschitz continuity and holds, for example, in $\ell_p$ regression with $p\in(1,2)$ \cite{AASS26}. It was studied in universal methods for convex problems \cite{Nes15} and in regularization methods \cite{CGT17}. For block methods, it was considered with random block selection \cite{MG25}, in a cyclic, partially derivative-free method for nonseparable composite problems \cite{Ama25} and, closest to the present work, by Amaral, Andreani, Secchin and Silva \cite{AASS26}, whose method combines free choice of the blocks, quadratic regularization and a stopping test that certifies approximate stationarity.

Block decomposition has received little attention in multiobjective optimization. For convex problems with a compact decision space, Soriano, Mitkovski and Wiecek \cite{SWM25} restricted the problem to lines, in the spirit of coordinate descent, and proved that the resulting sets converge in the Painlevé--Kuratowski sense; their approach is not a descent method. Tran and Vicente \cite{TV26} proposed a stochastic method that alternates between objectives and blocks of variables and obtained convergence rates in the convex and nonconvex cases, but its steps do not follow a descent direction common to all the objectives. To the best of our knowledge, there is no deterministic block method for nonconvex multiobjective optimization that uses common descent directions, has worst-case complexity bounds and stops only at approximately Pareto stationary points.

The difficulty already appears with two objectives and two variables. Consider $f_1(x)=x_1+\frac{1}{2}\|x\|^2$ and $f_2(x)=x_2+\frac{1}{2}\|x\|^2$, with blocks $\{x_1\}$ and $\{x_2\}$. At $x=0$, each block is stationary, since, in each coordinate, some convex combination of the partial gradients vanishes. Nevertheless, $x=0$ is not Pareto stationary, and every step that changes only one coordinate increases at least one of the objectives. This cannot happen in the differentiable scalar case, in which the gradient vanishes exactly when all the partial gradients vanish; this is the property that allows Amaral et al.\ \cite{AASS26} to use a stopping test checked block by block. With two or more objectives, Pareto stationarity requires a single convex combination of the gradients that vanishes in all the variables, whereas different blocks may be stationary with different combinations.

This phenomenon is the multiobjective counterpart of a known limitation of coordinate descent for nonsmooth and nonseparable functions \cite{War63,Tse01}. Indeed, when all the objectives have the same value at $x$, Pareto stationarity coincides with Clarke stationarity of the function $\max_if_i$ \cite[Prop.~2.3.12]{Cla83}. In the example above, the function $\max\{x_1,x_2\}+\frac{1}{2}\|x\|^2$ is minimized at $0$ along each coordinate, although $0$ is not stationary.

We call proper blocks the blocks of a partition of the variables, as opposed to the block formed by all of them. We show that this difficulty is not an isolated case. For any number of variables, there are problems with two objectives and an open set of points that are not Pareto stationary, from which no sequence of proper-block steps that does not increase any objective reduces the Pareto stationarity measure (Proposition~\ref{prop:impossibilidade}). Consequently, no method that uses only proper blocks and is monotone in all the objectives, in the sense that $f_i(x^{k+1})\leq f_i(x^k)$ for all $i$ and $k$, reaches Pareto stationarity in general, whatever the choice of directions and blocks. Extending block descent to multiobjective optimization therefore requires an ingredient that is unnecessary in the differentiable scalar case: a step outside the proper blocks or a nonmonotone acceptance rule.

In this paper, we propose a block coordinate descent method for multiobjective optimization under H\"older conditions that includes such an ingredient. At each iteration, one block, chosen freely, is updated along a common descent direction obtained from a regularized min-max model of the objectives restricted to that block. When the stationarity measures of all the proper blocks do not exceed the tolerance $\epsilon$ and the iterate is not yet $\epsilon$-Pareto stationary, the method takes a step in all the variables. Among the ways of dealing with the impossibility of Proposition~\ref{prop:impossibilidade}, this is the natural one for a monotone method, since the stopping test already requires the full gradient. Since termination is decided by the measure computed in all the variables, the method stops only at $\epsilon$-Pareto stationary points. The contributions of this work are the following.
\begin{itemize}
\item We show that no method restricted to proper blocks and monotone in all the objectives reaches Pareto stationarity in general, even with two objectives and any number of variables, and that the blockwise measures fail as a stopping test on an open set (Proposition~\ref{prop:impossibilidade} and Remark~\ref{obs:impossib-geral}). Since the objectives of the example are smooth and satisfy the assumptions of \cite{AASS26}, the failure cannot be removed by regularity assumptions and requires a change in the algorithm.
\item We propose a step in all the variables, taken only when no proper block exceeds the tolerance, with which the method stops only at $\epsilon$-Pareto stationary points. The bound of $O(\epsilon^{-(\beta+1)/\beta})$ iterations is of the same order as the bound of \cite{AASS26} for the scalar case and as the bound of \cite{CDM20} for methods that update all the variables. Since the step in all the variables reuses the partial gradients already computed at the current iterate, the number of gradient equivalents is at most the iteration bound plus one, so that, in the worst case, the method uses no more gradients than a method that updates all the variables with the same number of iterations (Theorem~\ref{teo:avaliacoes}(e)).
\item We analyze a family of update rules for the regularization parameter that includes the single-parameter rule of \cite{AASS26} and a new rule, with one parameter per block, for which we prove evaluation bounds (Theorem~\ref{teo:avaliacoes}(b) and (c)) and per-block bounds (Remark~\ref{obs:porbloco}).
\end{itemize}

The structure of the block step and the outline of the descent and complexity proofs come from \cite{AASS26}, the block subproblem is a regularized version of the direction of \cite{FS00}, and the step in all the variables corresponds to the method of \cite{CDM20} with $p=1$. The differences with respect to these works are pointed out after each result.

The rest of the paper is organized as follows. Section~\ref{sec:preliminares} introduces the notation, the blockwise stationarity measures and the assumptions. Section~\ref{sec:metodo} studies the block subproblem and contains the impossibility result and Algorithm~\ref{alg:bcd}. Section~\ref{sec:convergencia} shows that each iteration terminates after finitely many trials and establishes convergence and the complexity bounds. Section~\ref{sec:conclusoes} presents the conclusions, and Appendix~\ref{ap:provas} contains the proofs adapted from the scalar case.

\section{Preliminaries}\label{sec:preliminares} This section introduces the notation, the stationarity measures used by the method, and the assumptions of the analysis.

\subsection{Notation and problem formulation}
\label{sec:notacao}
Throughout the text, $\langle\cdot,\cdot\rangle$ and $\|\cdot\|$ denote the Euclidean inner product and norm, respectively, $\|\cdot\|_\infty$ denotes the maximum norm, and
\[
  \Delta^{r-1}:=\Big\{\gamma\in\mathbb{R}^r_+ :\ \sum_{i=1}^r\gamma_i=1\Big\}
\]
denotes the unit simplex. For $v,w\in\mathbb{R}^n$, the inequalities $v\leq w$ and $v<w$ are understood componentwise, and $e_i$ denotes the $i$-th canonical vector. For a finite set $S\subset\mathbb{R}^n$, $\operatorname{conv}S$ denotes its convex hull. For a symmetric matrix $A$, $\lambda_{\min}(A)$ denotes its smallest eigenvalue, and $\|A\|$ denotes the matrix norm induced by the Euclidean norm. 

We consider the problem
\begin{equation}\label{eq:MOP}
  \min_{x\in\mathbb{R}^n}\ F(x)=(f_1(x),\dots,f_r(x)),
\end{equation}
understood in the Pareto sense, that is, with respect to the partial order induced by $\mathbb{R}^r_+$. Each $f_i:\mathbb{R}^n\to\mathbb{R}$ is continuously differentiable and possibly nonconvex. We assume $r\geq 1$. The case of interest is $r\geq 2$; the case $r=1$ is used only for comparison with the scalar case.

The variables are partitioned into $q\geq 2$ blocks, called proper blocks, of dimensions $n_1,\dots,n_q\geq 1$ with $n_1+\cdots+n_q=n$. The partition is described by matrices $U_j\in\mathbb{R}^{n\times n_j}$, $j=1,\dots,q$, such that $[U_1\ \cdots\ U_q]$ is a permutation matrix. In particular,
\[
  U_j^TU_j=I_{n_j},\qquad U_j^TU_l=0\ \ (j\neq l),\qquad \sum_{j=1}^q U_jU_j^T=I_n .
\]
We denote by $J_j\subset\{1,\dots,n\}$ the set of indices of the columns of $I_n$ that form $U_j$, so that $v_{(j)}$ consists of the components $v_l$, $l\in J_j$.
We also introduce the block $j=0$, which represents the whole space, with $n_0:=n$ and $U_0:=I_n$. Thus, unless otherwise stated, $j$ ranges over the set $\{0,1,\dots,q\}$ throughout the text. For $v\in\mathbb{R}^n$, we write
\[
  v_{(j)}:=U_j^Tv\qquad\text{and}\qquad\nabla_{(j)}f_i(x):=U_j^T\nabla f_i(x),
\]
so that $v_{(0)}=v$ and $\nabla_{(0)}f_i(x)=\nabla f_i(x)$. We have $\|U_js\|=\|s\|$ for all $s\in\mathbb{R}^{n_j}$ and, since the proper blocks partition the coordinates,
\begin{equation}\label{eq:supnorm-blocos}
  \|v\|_\infty=\max_{1 \leq j \leq q}\|v_{(j)}\|_\infty,\qquad\forall v\in\mathbb{R}^n .
\end{equation}

\begin{definition}\label{def:pareto-estacionario}
We say that $x^*\in\mathbb{R}^n$ is \emph{Pareto stationary} for \eqref{eq:MOP} if
\[
  0\in\operatorname{conv}\{\nabla f_1(x^*),\dots,\nabla f_r(x^*)\},
\]
that is, if there exists $\gamma\in\Delta^{r-1}$ such that $\sum_{i=1}^r\gamma_i\nabla f_i(x^*)=0$.
\end{definition}

Equivalently, $x^*$ is Pareto stationary if and only if there is no $d\in\mathbb{R}^n$ such that
\[
  \langle\nabla f_i(x^*),d\rangle<0,\qquad i=1,\dots,r,
\]
that is, if there is no common descent direction for all the objectives, which is the notion of Pareto criticality introduced in \cite{FS00}; see \cite[Lemma~3.6]{Sam26}. We say that $x^*$ is a \emph{local weak Pareto minimizer} of \eqref{eq:MOP} if there is a neighborhood $V$ of $x^*$ such that no $y\in V$ satisfies $f_i(y)<f_i(x^*)$ for all $i=1,\dots,r$, and a \emph{weak Pareto minimizer} if this holds with $V=\mathbb{R}^n$. Pareto stationarity is a necessary condition for $x^*$ to be a local weak Pareto minimizer, and it is also sufficient for $x^*$ to be a weak Pareto minimizer when all the $f_i$ are convex \cite[Theorem~3.1]{FGS09}.

\subsection{Stationarity measures}
The following measures are used by Algorithm~\ref{alg:bcd} to decide which blocks are skipped and when to stop.

\begin{definition}\label{def:Gj}
Given $x\in\mathbb{R}^n$ and $j\in\{0,1,\dots,q\}$, we define the \emph{Pareto stationarity measure on block $j$} as
\begin{equation}\label{eq:Gj}
  G_j(x):=\min_{\gamma\in\Delta^{r-1}}\Big\|\sum_{i=1}^r\gamma_i\nabla_{(j)}f_i(x)\Big\|_\infty .
\end{equation}
For $j=0$, the block consists of all the variables, and we write
\[
  G(x):=G_0(x)=\min_{\gamma\in\Delta^{r-1}}\Big\|\sum_{i=1}^r\gamma_i\nabla f_i(x)\Big\|_\infty ,
\]
which we call the \emph{Pareto stationarity measure}.
\end{definition}

The minimum in \eqref{eq:Gj} is attained, since the function $\gamma\mapsto\|\sum_i\gamma_i\nabla_{(j)}f_i(x)\|_\infty$ is continuous and $\Delta^{r-1}$ is compact. In particular, $G(x)=0$ if and only if $x$ is Pareto stationary. For $j\geq1$, the equality $G_j(x)=0$ means that some convex combination of the partial gradients with respect to the variables of block $j$ vanishes. This combination may depend on $j$, and the consequences of this are discussed in Remark~\ref{obs:medidas-bloco}.

The measure $G$ is the analogue, in the maximum norm, of the Pareto stationarity measure of \cite{Sam26}, which is defined with the Euclidean norm. The choice of norm is discussed in Remark~\ref{obs:norma}.

\begin{definition}\label{def:eps-pareto}
Given $\epsilon\geq0$, we say that $x\in\mathbb{R}^n$ is \emph{$\epsilon$-Pareto stationary} for \eqref{eq:MOP} if $G(x)\leq\epsilon$.
\end{definition}

The following lemma collects the properties of the measures $G_j$ used in the sequel. Item~(a) shows that, if the measure of some proper block exceeds the tolerance, then the point is not yet $\epsilon$-Pareto stationary, and item~(b) shows that, with a single objective, the blockwise measures determine $G$. Item~(c) is used in the convergence analysis (Corollary~\ref{cor:convergencia}), and item~(d) relates $G$ to the measure defined with the Euclidean norm.

\begin{lemma}\label{lem:prop-Gj}
For all $x,y\in\mathbb{R}^n$, the following statements hold.
\begin{enumerate}
  \item[(a)] $G(x)\geq\max_{1 \leq j \leq q}G_j(x)$.
  \item[(b)] If $r=1$, equality holds in (a).
  \item[(c)] For every $j\in\{0,1,\dots,q\}$,
  \[
    |G_j(x)-G_j(y)| \leq \max_{1 \leq i \leq r}\|\nabla_{(j)}f_i(x)-\nabla_{(j)}f_i(y)\|_\infty .
  \]
  In particular, each $G_j$ is continuous.
  \item[(d)] $\displaystyle G(x)\ \leq\ \min_{\gamma\in\Delta^{r-1}}\Big\|\sum_{i=1}^r\gamma_i\nabla f_i(x)\Big\|\ \leq \ \sqrt n\,G(x)$.
\end{enumerate}
\end{lemma}

\begin{proof}\leavevmode
\begin{enumerate}
    \item[(a)] Let $\gamma^*$ be a minimizer in \eqref{eq:Gj} with $j=0$. For each $j\geq1$, the components of $v_{(j)}$ are components of $v$, so that $\|v_{(j)}\|_\infty\leq\|v\|_\infty$ for all $v\in\mathbb{R}^n$. Hence,
\[
  G_j(x)\leq\Big\|\sum_{i=1}^r\gamma_i^*\nabla_{(j)}f_i(x)\Big\|_\infty\leq\Big\|\sum_{i=1}^r\gamma_i^*\nabla f_i(x)\Big\|_\infty=G(x).
\]
\item[(b)] If $r=1$, then $\Delta^0=\{1\}$ and $G_j(x)=\|\nabla_{(j)}f_1(x)\|_\infty$ for all $j$. By \eqref{eq:supnorm-blocos},
\[
  G(x)=\|\nabla f_1(x)\|_\infty=\max_{1 \leq j \leq q}\|\nabla_{(j)}f_1(x)\|_\infty=\max_{1 \leq j \leq q}G_j(x).
\]
\item[(c)] Let $M:=\max_{1 \leq i \leq r}\|\nabla_{(j)}f_i(x)-\nabla_{(j)}f_i(y)\|_\infty$. For each $\gamma\in\Delta^{r-1}$, the triangle inequality gives
\begin{align*}
\Big\|\sum_{i=1}^r\gamma_i\nabla_{(j)}f_i(x)\Big\|_\infty
&\leq\Big\|\sum_{i=1}^r\gamma_i\nabla_{(j)}f_i(y)\Big\|_\infty+\sum_{i=1}^r\gamma_i\|\nabla_{(j)}f_i(x)-\nabla_{(j)}f_i(y)\|_\infty\\
&\leq\Big\|\sum_{i=1}^r\gamma_i\nabla_{(j)}f_i(y)\Big\|_\infty+M .
\end{align*}
Taking $\gamma$ as a minimizer in \eqref{eq:Gj} at the point $y$, we obtain $G_j(x)\leq G_j(y)+M$. Interchanging the roles of $x$ and $y$, we obtain $G_j(y)\leq G_j(x)+M$, which proves the inequality. The continuity of $G_j$ follows from the continuity of the gradients $\nabla f_i$.
\item[(d)] For every $\gamma\in\Delta^{r-1}$,
\[
  \Big\|\sum_{i=1}^r\gamma_i\nabla f_i(x)\Big\|_\infty \leq \Big\|\sum_{i=1}^r\gamma_i\nabla f_i(x)\Big\|\leq\sqrt n\,\Big\|\sum_{i=1}^r\gamma_i\nabla f_i(x)\Big\|_\infty .
\]
Since these inequalities hold for every $\gamma$, they are preserved when the minimum over $\gamma$ is taken in each term, which yields item~(d).
\end{enumerate}
\end{proof}

\begin{remark}\label{obs:medidas-bloco}
For $r\geq2$, the inequality in Lemma~\ref{lem:prop-Gj}(a) may be strict, and it may happen that
\[
  G_j(x)=0,\quad j=1,\dots,q,\qquad\text{with}\qquad G(x)>0,
\]
as in the example of the Introduction. Therefore, the measures $G_1,\dots,G_q$ are not sufficient to ensure Pareto stationarity. The continuity of the optimal value of the subproblem defining the multiobjective steepest descent direction, analogous to Lemma~\ref{lem:prop-Gj}(c), is established in \cite[Lemma~1]{FS00}.
\end{remark}

\begin{remark}\label{obs:norma}
We adopt the maximum norm in \eqref{eq:Gj} for compatibility with the scalar case. By Lemma~\ref{lem:prop-Gj}(b), when $r=1$ we have
\[
  G(x)=\|\nabla f_1(x)\|_\infty=\max_{1 \leq j \leq q}G_j(x),
\]
so that the measure $G_j$ and the stopping criterion coincide with those of \cite{AASS26}. With the Euclidean norm, we would only have
\[
  G(x)\leq\sqrt q\,\max_{1 \leq j \leq q}G_j(x),
\]
and Algorithm~\ref{alg:bcd} could resort to block~0 even when $r=1$.

Moreover, with the maximum norm, the computation of $G_j$ reduces to a linear programming problem. For fixed vectors $v_1,\dots,v_r\in\mathbb{R}^{n_j}$, the value $\min_{\gamma\in\Delta^{r-1}}\|\sum_i\gamma_iv_i\|_\infty$ is the optimal value of the problem with $r+1$ variables
\[
  \min_{t\in\mathbb{R},\,\gamma\in\Delta^{r-1}}\ t\quad\text{s.t.}\quad -t\leq\sum_{i=1}^r\gamma_i\,[v_i]_l \leq t,\qquad l=1,\dots,n_j,
\]
where $[v_i]_l$ denotes the $l$-th component of $v_i$.

If the maximum norm in \eqref{eq:Gj} is replaced by the Euclidean norm, all the results of Sections~\ref{sec:algoritmo} to~\ref{sec:complexidade} remain valid with the same constants. Indeed, the proofs only use the fact that the norm of a subvector does not exceed that of the vector and, in Lemma~\ref{lem:cota-passo}, the inequality $\|\cdot\|_\infty\leq\|\cdot\|$, which becomes an equality. In this case, the stopping criterion directly controls the Pareto stationarity measure of \cite{Sam26}, without the factor $\sqrt n$ of Lemma~\ref{lem:prop-Gj}(d), and the only loss is the coincidence with \cite{AASS26} when $r=1$.
\end{remark}
\subsection{Assumptions}
\label{sec:hipoteses}

This subsection collects the assumptions used in the analysis. Assumption~\ref{hip:holder} bounds the error of the linear approximation of each objective along each block, Assumption~\ref{hip:matrizes} bounds the matrices of the quadratic models, and Assumption~\ref{hip:limitacao} requires one of the objectives to be bounded below on its initial level set.

\begin{assumption}\label{hip:holder}
For each $i\in\{1,\dots,r\}$ and $j\in\{0,1,\dots,q\}$, there exist $L_{i,j}>0$ and $\beta_{i,j}\in(0,1]$ such that
\[
  f_i(x+U_js_{(j)}) \leq f_i(x)+\langle\nabla_{(j)}f_i(x),s_{(j)}\rangle+L_{i,j}\|s_{(j)}\|^{1+\beta_{i,j}},
  \qquad\forall\,x\in\mathbb{R}^n,\ \forall\,s_{(j)}\in\mathbb{R}^{n_j}.
\]
\end{assumption}

\begin{remark}\label{obs:holder}
Assumption~\ref{hip:holder} is implied by the H\"older continuity of the partial gradients. Indeed, suppose that there exists $H_{i,j}>0$ such that
\[
  \|\nabla_{(j)}f_i(x+U_js_{(j)})-\nabla_{(j)}f_i(x)\| \leq H_{i,j}\|s_{(j)}\|^{\beta_{i,j}},
  \qquad\forall\,x\in\mathbb{R}^n,\ \forall\,s_{(j)}\in\mathbb{R}^{n_j}.
\]
For $g(t):=f_i(x+tU_js_{(j)})$, we have
\[
  g'(t)=\langle\nabla_{(j)}f_i(x+tU_js_{(j)}),s_{(j)}\rangle
\]
and, by the Cauchy--Schwarz inequality,
\begin{align*}
  &f_i(x+U_js_{(j)})-f_i(x)-\langle\nabla_{(j)}f_i(x),s_{(j)}\rangle\\
  &\qquad=\int_0^1\langle\nabla_{(j)}f_i(x+tU_js_{(j)})-\nabla_{(j)}f_i(x),s_{(j)}\rangle\,dt\\
  &\qquad\leq H_{i,j}\|s_{(j)}\|^{1+\beta_{i,j}}\int_0^1t^{\beta_{i,j}}\,dt\\
  &\qquad=\frac{H_{i,j}}{1+\beta_{i,j}}\|s_{(j)}\|^{1+\beta_{i,j}} .
\end{align*}
Hence, Assumption~\ref{hip:holder} holds with $L_{i,j}=H_{i,j}/(1+\beta_{i,j})$, as in \cite{AASS26}. The converse is false, since every concave function of class $C^1$ satisfies the inequality for any $L_{i,j}>0$, without its gradient being necessarily H\"older continuous.

Since $\|U_js_{(j)}\|=\|s_{(j)}\|$ and $\langle\nabla f_i(x),U_js_{(j)}\rangle=\langle\nabla_{(j)}f_i(x),s_{(j)}\rangle$, Assumption~\ref{hip:holder} with $j=0$ implies the same inequality for every $j\in\{1,\dots,q\}$, with $L_{i,j}=L_{i,0}$ and $\beta_{i,j}=\beta_{i,0}$. The blockwise constants only allow sharper values.

In the scalar case, \cite{AASS26} requires the condition only on the proper blocks. Here, the condition with $j=0$ is used only in the steps on block~0, that is, in the steps outside the proper blocks, which a monotone method generally requires when $r\geq2$, as shown in Section~\ref{sec:algoritmo}.
\end{remark}
\begin{assumption}\label{hip:matrizes}
There exists $B\geq1$ such that, for all $i\in\{1,\dots,r\}$, $j\in\{0,1,\dots,q\}$ and $x\in\mathbb{R}^n$, the matrix $B_{i,(j)}(x)\in\mathbb{R}^{n_j\times n_j}$ is symmetric and $\|B_{i,(j)}(x)\|\leq B$, where $\|\cdot\|$ denotes the matrix norm induced by the Euclidean norm.
\end{assumption}

The matrices $B_{i,(j)}(x)$ define the quadratic term of the model of $f_i$ on block $j$, given in \eqref{eq:modelo}. Examples of matrices that satisfy Assumption~\ref{hip:matrizes} are the zero matrix, bounded quasi-Newton approximations, and the partial Hessian $U_j^T\nabla^2f_i(x)U_j$ when $f_i$ is of class $C^2$ with bounded Hessian.

\begin{assumption}\label{hip:limitacao}
There exist $i_0\in\{1,\dots,r\}$ and $f_{\mathrm{low}}\in\mathbb{R}$ such that $f_{i_0}(x)\geq f_{\mathrm{low}}$ for all $x$ in the level set
\[
  \mathcal{L}_0:=\{x\in\mathbb{R}^n:\ f_{i_0}(x)\leq f_{i_0}(x^0)\},
\]
where $x^0$ is the initial point of Algorithm~\ref{alg:bcd}.
\end{assumption}

Assumption~\ref{hip:limitacao} holds whenever some $f_i$ is bounded below on $\mathbb{R}^n$.

\section{The method}\label{sec:metodo}
This section describes the subproblem solved on each block and presents the algorithm.
\subsection{The block subproblem}
\label{sec:subproblema}

In this subsection, we describe the subproblem solved by Algorithm~\ref{alg:bcd} whenever a step is computed on a block $j\in\{0,1,\dots,q\}$. The idea is to approximate the variation of each objective $f_i$ along block $j$ by a quadratic model, add a regularization term that controls the step size, and seek a step that decreases all these models simultaneously, which leads to a min-max problem.

Given the current iterate $x^k$, the block $j$ and the regularization parameter $\sigma\geq1$ provided by Algorithm~\ref{alg:bcd}, we define, for $z\in\mathbb{R}^{n_j}$ and $i=1,\dots,r$,
\begin{equation}\label{eq:modelo}
\varphi_i(z):=\langle\nabla_{(j)}f_i(x^k),z\rangle+\frac{1}{2}\langle B_{i,(j)}(x^k)z,z\rangle+\sigma\|z\|^2,
\end{equation}
which is a regularized model of $f_i(x^k+U_jz)-f_i(x^k)$. To simplify the notation, we omit the dependence of $\varphi_i$ on $k$, $j$ and $\sigma$ and, with $j$ fixed, we write the step $s_{(j)}\in\mathbb{R}^{n_j}$ simply as $s$. The step on block $j$ is an approximate solution, in the sense of Definition~\ref{def:aceitavel}, of the problem
\begin{equation}\label{eq:subproblema}
  \min_{z\in\mathbb{R}^{n_j}}\ \max_{1 \leq i \leq r}\ \varphi_i(z).
\end{equation}
For $r=1$, problem \eqref{eq:subproblema} is the subproblem of \cite{AASS26}. For $j=0$ and $B_{i,(0)}\equiv0$, it is the subproblem defining the steepest descent direction of Fliege and Svaiter \cite{FS00}, up to the scaling factor $1/(2\sigma)$.

It is not necessary to solve \eqref{eq:subproblema} exactly. It suffices to find a step that does not increase any of the models and is approximately stationary, in the sense of the following definition.

\begin{definition}\label{def:aceitavel}
Given $\theta>0$, we say that $s_{(j)}\in\mathbb{R}^{n_j}$ and $\gamma\in\Delta^{r-1}$ form an \emph{acceptable pair} if
\begin{align}
  \varphi_i(s_{(j)})& \leq 0,\qquad i=1,\dots,r, \label{eq:aceit1}\\
\Big\|\sum_{i=1}^r\gamma_i\nabla\varphi_i(s_{(j)})\Big\|&\leq\theta\|s_{(j)}\|. \label{eq:aceit2}
\end{align}
\end{definition}

For $r=1$, conditions \eqref{eq:aceit1} and \eqref{eq:aceit2} reduce to conditions (2) and (3) of \cite{AASS26}. For $j=0$, they correspond to conditions (5)--(7) of \cite{CDM20}, without constraints, with $p=1$ and the model
$$
  M_i(x^k+z)=f_i(x^k)+\langle\nabla f_i(x^k),z\rangle+\frac{1}{2}\langle B_{i,(0)}(x^k)z,z\rangle .
$$

The following lemma shows that, when the matrices $B_{i,(j)}(x^k)+2\sigma I_{n_j}$ are positive definite, problem \eqref{eq:subproblema} has a unique solution, characterized by the vanishing of a convex combination of the gradients of the models, and that this solution yields an acceptable pair. Since Step~5 of Algorithm~\ref{alg:bcd} ensures this positive definiteness, the computation of the step in Algorithm~\ref{alg:bcd} is always well defined.

\begin{lemma}\label{lem:subproblema}
Suppose that $H_i:=B_{i,(j)}(x^k)+2\sigma I_{n_j}$ is symmetric and positive definite for all $i=1,\dots,r$, and let $\phi(z):=\max_{1 \leq i \leq r}\varphi_i(z)$. Then the following statements hold.
\begin{enumerate}
  \item[(a)] The function $\phi$ is strongly convex and has a unique minimizer $\bar s\in\mathbb{R}^{n_j}$.
  \item[(b)] There exists $\gamma^k\in\Delta^{r-1}$, with $\gamma^k_i=0$ whenever $\varphi_i(\bar s)<\phi(\bar s)$, such that
  \begin{equation}\label{eq:kkt-sub}
    \sum_{i=1}^r\gamma_i^k\nabla\varphi_i(\bar s)=0.
  \end{equation}
  \item[(c)] $\varphi_i(\bar s)\leq0$ for all $i=1,\dots,r$.
  \item[(d)] For every $\theta>0$, $(\bar s,\gamma^k)$ is an acceptable pair in the sense of Definition~\ref{def:aceitavel}.
\end{enumerate}
\end{lemma}
\begin{proof}
Let $g_i:=\nabla_{(j)}f_i(x^k)$, so that
$$
  \varphi_i(z)=\langle g_i,z\rangle+\frac{1}{2}\langle H_iz,z\rangle\qquad\text{and}\qquad\nabla\varphi_i(z)=g_i+H_iz,
$$
and let $\mu_i:=\lambda_{\min}(H_i)$, which is positive because $H_i$ is symmetric and positive definite, and $\mu:=\min_i\mu_i>0$.
\begin{enumerate}
\item[(a)] Given $z_1,z_2\in\mathbb{R}^{n_j}$, $\lambda\in[0,1]$ and $\bar z:=\lambda z_1+(1-\lambda)z_2$, since $\varphi_i$ is quadratic,
$$
  \varphi_i(\bar z)=\lambda\varphi_i(z_1)+(1-\lambda)\varphi_i(z_2)-\tfrac{1}{2}\lambda(1-\lambda)\langle H_i(z_1-z_2),z_1-z_2\rangle .
$$
Using $\varphi_i(z_1)\leq\phi(z_1)$, $\varphi_i(z_2)\leq\phi(z_2)$ and $\langle H_i(z_1-z_2),z_1-z_2\rangle\geq\mu\|z_1-z_2\|^2$, and taking the maximum over $i$, we obtain
$$
  \phi(\bar z)\leq\lambda\phi(z_1)+(1-\lambda)\phi(z_2)-\tfrac{\mu}{2}\lambda(1-\lambda)\|z_1-z_2\|^2 ,
$$
that is, $\phi$ is $\mu$-strongly convex.

The function $\phi$ is continuous and satisfies $\phi\geq\varphi_1$, with
$$
  \varphi_1(z)\geq-\|g_1\|\|z\|+\frac{\mu_1}{2}\|z\|^2\to\infty\qquad\text{as }\|z\|\to\infty .
$$
Hence, $\phi$ is coercive and, by the Weierstrass theorem, has a global minimizer $\bar s$.

For uniqueness, suppose that there exists another minimizer $\bar s'\neq\bar s$. Since $\phi(\bar s')=\phi(\bar s)$, strong convexity with $z_1=\bar s$, $z_2=\bar s'$ and $\lambda=\frac{1}{2}$ gives
$$
  \phi\big(\tfrac{1}{2}(\bar s+\bar s')\big)\leq\phi(\bar s)-\tfrac{\mu}{8}\|\bar s-\bar s'\|^2<\phi(\bar s),
$$
which contradicts the minimality of $\bar s$.

\item[(b)] Let $I:=\{i:\varphi_i(\bar s)=\phi(\bar s)\}$ and $d\in\mathbb{R}^{n_j}$. For $i\notin I$, we have $\varphi_i(\bar s)<\phi(\bar s)$. Since the functions $t\mapsto\varphi_i(\bar s+td)$ are continuous, there exists $t_0>0$ such that, for $t\in(0,t_0)$, the maximum defining $\phi(\bar s+td)$ is attained in $I$. Since $\varphi_i(\bar s)=\phi(\bar s)$ for $i\in I$,
$$
  \frac{\phi(\bar s+td)-\phi(\bar s)}{t}=\max_{i\in I}\Big[\langle\nabla\varphi_i(\bar s),d\rangle+\frac{t}{2}\langle H_id,d\rangle\Big],\qquad\forall\,t\in(0,t_0).
$$
Since $\bar s$ minimizes $\phi$, the left-hand side is nonnegative. Letting $t\downarrow0$, we obtain
\begin{equation}\label{eq:dir-deriv}
  \max_{i\in I}\langle\nabla\varphi_i(\bar s),d\rangle\geq0,\qquad\forall\,d\in\mathbb{R}^{n_j}.
\end{equation}

Let $K:=\operatorname{conv}\{\nabla\varphi_i(\bar s):i\in I\}$, which is compact and convex. Suppose, by contradiction, that $0\notin K$, and let $\bar v\neq0$ be the projection of $0$ onto $K$. By the characterization of the projection,
$$
  \langle v,\bar v\rangle\geq\|\bar v\|^2,\qquad\forall\,v\in K.
$$
Taking $d:=-\bar v$, we obtain
$$
  \max_{i\in I}\langle\nabla\varphi_i(\bar s),d\rangle\leq-\|\bar v\|^2<0,
$$
which contradicts \eqref{eq:dir-deriv}. Hence, $0\in K$, that is, there exist $\gamma_i^k\geq0$, $i\in I$, with $\sum_{i\in I}\gamma_i^k=1$ and $\sum_{i\in I}\gamma_i^k\nabla\varphi_i(\bar s)=0$. Setting $\gamma_i^k:=0$ for $i\notin I$, we obtain \eqref{eq:kkt-sub}.

\item[(c)] For each $i\in\{1,\dots,r\}$, by the definition of $\phi$ and the minimality of $\bar s$,
$$
  \varphi_i(\bar s)\leq\phi(\bar s)\leq\phi(0)=0,
$$
since $\varphi_l(0)=0$ for all $l$.

\item[(d)] Let $\theta>0$. By item~(b), $\gamma^k\in\Delta^{r-1}$, and condition \eqref{eq:aceit1} is item~(c). As for \eqref{eq:aceit2}, by \eqref{eq:kkt-sub},
$$
  \Big\|\sum_{i=1}^r\gamma_i^k\nabla\varphi_i(\bar s)\Big\|=0\leq\theta\|\bar s\|.
$$
Hence, $(\bar s,\gamma^k)$ is an acceptable pair.
\end{enumerate}
\end{proof}

\begin{remark}\label{obs:calculo-passo}
By Lemma~\ref{lem:subproblema}(d), the exact minimizer of \eqref{eq:subproblema} always yields an acceptable pair. Its computation depends on the matrices $B_{i,(j)}(x^k)$.
\begin{enumerate}
  \item[(i)] If $B_{i,(j)}(x^k)=B_{(j)}(x^k)$ for all $i$, with $H:=B_{(j)}(x^k)+2\sigma I_{n_j}$ positive definite, the minimizer has the form
  $$
    \bar s=-H^{-1}\sum_{i=1}^r\gamma_i^k\nabla_{(j)}f_i(x^k),
  $$
  where $\gamma^k$ is a solution of the $r$-dimensional quadratic problem
  \begin{equation}\label{eq:dual-sub}
    \min_{\gamma\in\Delta^{r-1}}\ \frac{1}{2}\Big\langle H^{-1}\sum_{i=1}^r\gamma_i\nabla_{(j)}f_i(x^k),\ \sum_{i=1}^r\gamma_i\nabla_{(j)}f_i(x^k)\Big\rangle .
  \end{equation}
  Indeed, since
  $$
    \max_{1 \leq i \leq r}\varphi_i(z)=\max_{\gamma\in\Delta^{r-1}}\sum_{i=1}^r\gamma_i\varphi_i(z),
  $$
  and the function $(z,\gamma)\mapsto\sum_i\gamma_i\varphi_i(z)$ is convex in $z$ and linear in $\gamma$, with $\Delta^{r-1}$ compact, Sion's minimax theorem allows us to interchange the minimum over $z$ and the maximum over $\gamma$. For fixed $\gamma$, the minimum of $\sum_i\gamma_i\varphi_i$ is attained at
  $$
    z(\gamma)=-H^{-1}\sum_{i=1}^r\gamma_i\nabla_{(j)}f_i(x^k),
  $$
  with value
  $$
    -\frac{1}{2}\Big\langle H^{-1}\sum_{i=1}^r\gamma_i\nabla_{(j)}f_i(x^k),\ \sum_{i=1}^r\gamma_i\nabla_{(j)}f_i(x^k)\Big\rangle ,
  $$
  and maximizing this value over $\gamma$ is equivalent to solving \eqref{eq:dual-sub}. If $\gamma^k$ is a solution of \eqref{eq:dual-sub}, then $(\bar s,\gamma^k)$ is a saddle point. Hence, $\bar s$ minimizes $\sum_i\gamma_i^k\varphi_i$, which is strongly convex, and therefore $\bar s=z(\gamma^k)$.
  \item[(ii)] If the matrices $B_{i,(j)}(x^k)$ depend on $i$, subproblem \eqref{eq:subproblema} is equivalent to the convex quadratically constrained problem
  $$
    \min_{t\in\mathbb{R},\,z\in\mathbb{R}^{n_j}}\ t\quad\text{s.t.}\quad\varphi_i(z)\leq t,\qquad i=1,\dots,r,
  $$
  which can be solved by standard convex optimization methods.
\end{enumerate}
\end{remark}

\subsection{The algorithm}
\label{sec:algoritmo}

In this subsection, we present Algorithm~\ref{alg:bcd}. First, we show why, with two or more objectives, a block method that does not increase any objective needs steps outside the proper blocks. In the scalar case, the condition $G_j(x^k)\leq\epsilon$ for all proper blocks ensures that $x^k$ is $\epsilon$-stationary (Lemma~\ref{lem:prop-Gj}(b)), and the method of \cite{AASS26} can stop using only these blocks. For $r\geq2$, this no longer holds, even with $\epsilon=0$. The following proposition shows that, for any number of variables $n\geq2$, there exists a problem with an open region of points that are not Pareto stationary, from which no sequence of steps on proper blocks that do not increase any objective reduces $G$.

\begin{proposition}\label{prop:impossibilidade}
Let $n\geq r\geq2$, $x=(z,y)\in\mathbb{R}^r\times\mathbb{R}^{n-r}$, $\psi(y):=\frac{1}{2}\max\{\|y\|-1,0\}^2$ and
$$
  f_i(x):=z_i+\tfrac{1}{2}\|z\|^2+\psi(y),\qquad i=1,\dots,r,
$$
where $y$ and $\psi$ do not appear when $n=r$. Consider a partition of the variables into $q\geq2$ proper blocks such that none of the index sets $J_1,\dots,J_q$ contains $\{1,\dots,r\}$, and let
$$
  \bar\Omega:=\Big\{x\in\mathbb{R}^n:\ z\leq0,\ \sum_{i\in J_j,\,i\leq r}|z_i|\leq1,\ j=1,\dots,q,\ \|y\|\leq1\Big\}
$$
and
$$
  \Omega:=\Big\{x\in\mathbb{R}^n:\ z<0,\ \sum_{i\in J_j,\,i\leq r}|z_i|<1,\ j=1,\dots,q,\ \|y\|<1,\ \sum_{i=1}^rz_i\neq-1\Big\}.
$$
Then the following statements hold.
\begin{enumerate}
  \item[(a)] If $x\in\bar\Omega$, $j\in\{1,\dots,q\}$, $s\in\mathbb{R}^{n_j}$ and $f_i(x+U_js)\leq f_i(x)$ for all $i$, then $x+U_js\in\bar\Omega$ and $x+U_js$ has the same $z$-component as $x$.
  \item[(b)] $G_j(x)=0$ for all $x\in\bar\Omega$ and $j=1,\dots,q$.
  \item[(c)] A point $x\in\bar\Omega$ is Pareto stationary if and only if $\sum_{i=1}^rz_i=-1$. Moreover, $\Omega\subset\bar\Omega$ is open, nonempty, and contains no Pareto stationary points.
  \item[(d)] If $x^0\in\Omega$ and $x^{k+1}=x^k+U_{j_k}s^k$, with $j_k\in\{1,\dots,q\}$ and $f_i(x^{k+1})\leq f_i(x^k)$ for all $i$ and $k$, then $x^k\in\bar\Omega$ and $G(x^k)=G(x^0)>0$ for all $k$.
\end{enumerate}
\end{proposition}

\begin{proof}
The function $\psi$ is half the squared distance to the closed unit ball. Hence, $\psi\geq0$, and $\psi$ is of class $C^1$, with
$$
  \nabla\psi(y)=y-P(y),
$$
where $P$ is the projection onto this ball. In particular, $\psi(y)=0$ and $\nabla\psi(y)=0$ if and only if $\|y\|\leq1$. Therefore, for $x\in\bar\Omega$,
$$
  \nabla f_i(x)=(e_i+z,\,0),\qquad i=1,\dots,r,
$$
where $e_i$ denotes here the $i$-th canonical vector of $\mathbb{R}^r$.
\begin{enumerate}
  \item[(a)] Let $d:=U_js=(d_z,d_y)$, so that $d_l=0$ for $l\notin J_j$. Since $\psi(y)=0$, for $i=1,\dots,r$,
  $$
    f_i(x+d)-f_i(x)=d_i+\langle z,d_z\rangle+\tfrac{1}{2}\|d_z\|^2+\psi(y+d_y).
  $$
  Suppose that $f_i(x+d)\leq f_i(x)$ for all $i$ and, by contradiction, that $d_z\neq0$. Let
  $$
    a:=-\langle z,d_z\rangle-\tfrac{1}{2}\|d_z\|^2-\psi(y+d_y),
  $$
  so that $d_i\leq a$ for all $i\leq r$. Since $J_j$ does not contain $\{1,\dots,r\}$, there exists $i_0\leq r$ with $i_0\notin J_j$, and then $d_{i_0}=0$. Hence, $a\geq0$.

  Since $d_z\neq0$, the set $J_j\cap\{1,\dots,r\}$ is nonempty. Let $m:=\max_{i\in J_j,\,i\leq r}d_i$, which satisfies $m\leq a$, and $m_+:=\max\{m,0\}$. With $w:=-z$, we have $w\geq0$ and $\sum_{i\in J_j,\,i\leq r}w_i\leq1$, and therefore
  $$
    -\langle z,d_z\rangle=\sum_{i\in J_j,\,i\leq r}w_id_i\leq m_+ .
  $$
  Since $\psi\geq0$ and $d_z\neq0$, it follows that
  $$
    a\leq m_+-\tfrac{1}{2}\|d_z\|^2<m_+ .
  $$
  If $m\leq0$, then $a<0$, which contradicts $a\geq0$. If $m>0$, then $m\leq a<m$, which is also absurd. Hence, $d_z=0$.

  With $d_z=0$, we have $f_i(x+d)-f_i(x)=\psi(y+d_y)$ for all $i$, and the hypothesis gives $\psi(y+d_y)=0$, that is, $\|y+d_y\|\leq1$. Thus, $x+d\in\bar\Omega$ and has the same $z$-component as $x$.

  \item[(b)] Fix $j$ and $i_0\leq r$ with $i_0\notin J_j$, which exists by assumption. Define $\gamma_i:=-z_i$ for $i\in J_j$ with $i\leq r$, $\gamma_{i_0}:=1-\sum_{i\in J_j,\,i\leq r}|z_i|$, and $\gamma_i:=0$ for the remaining indices. Since $x\in\bar\Omega$, we have $\gamma\in\Delta^{r-1}$ and
  $$
    \sum_{i=1}^r\gamma_i\nabla f_i(x)=(\gamma+z,\,0),
  $$
  whose components in $J_j$ vanish. Hence, $G_j(x)=0$.

  \item[(c)] By Definition~\ref{def:pareto-estacionario}, $x\in\bar\Omega$ is Pareto stationary if and only if $\gamma+z=0$ for some $\gamma\in\Delta^{r-1}$. Since $z\leq0$, this holds if and only if $\sum_{i=1}^rz_i=-1$. The set $\Omega$ is open, as the intersection of finitely many open sets, and is contained in $\bar\Omega$. It is nonempty, since it contains $x=(-\delta(1,\dots,1),0)$ for every $\delta\in(0,1/r)$, and it contains no Pareto stationary points, since $\sum_{i=1}^rz_i\neq-1$ in $\Omega$.

  \item[(d)] Let $x^0=(z^0,y^0)\in\Omega$. If $x^k=(z^0,y^k)\in\bar\Omega$, item~(a), applied with $x=x^k$, $j=j_k$ and $s=s^k$, gives
  $$
    x^{k+1}=x^k+U_{j_k}s^k=(z^0,\,y^{k+1})\in\bar\Omega .
  $$
  By induction, $x^k=(z^0,y^k)\in\bar\Omega$ for all $k$. Since $\nabla f_i(x)=(e_i+z,\,0)$ on $\bar\Omega$,
  $$
    \nabla f_i(x^k)=(e_i+z^0,\,0)=\nabla f_i(x^0),\qquad i=1,\dots,r,
  $$
  and therefore
  $$
    G(x^k)=\min_{\gamma\in\Delta^{r-1}}\Big\|\sum_{i=1}^r\gamma_i\nabla f_i(x^k)\Big\|_\infty=G(x^0)>0,
  $$
  where the inequality follows from item~(c).
\end{enumerate}
\end{proof}
The functions of Proposition~\ref{prop:impossibilidade} satisfy Assumptions~\ref{hip:holder} and~\ref{hip:limitacao}. Indeed, $\nabla f_i(x)=(e_i+z,\,y-P(y))$ is Lipschitz continuous with constant $1$, since $I-P$ is nonexpansive, and Remark~\ref{obs:holder} gives Assumption~\ref{hip:holder} with
$$
  \beta_{i,j}=1,\qquad L_{i,j}=\tfrac{1}{2},\qquad i=1,\dots,r,\quad j=0,1,\dots,q.
$$
Moreover, since $\psi\geq0$,
$$
  f_i(x)\geq z_i+\tfrac{1}{2}\|z\|^2\geq z_i+\tfrac{1}{2}z_i^2\geq-\tfrac{1}{2},
$$
and Assumption~\ref{hip:limitacao} holds. Since $q\geq2$, the condition on $J_1,\dots,J_q$ does not restrict the partition, because it suffices to reindex the variables so that the first $r$ of them do not all belong to the same block. Therefore, for every partition into proper blocks, there exist a problem with $r=2$ that satisfies Assumptions~\ref{hip:holder} and~\ref{hip:limitacao} and an open set $\Omega$ of initial points from which no monotone method restricted to these blocks reaches an $\epsilon$-Pareto stationary point with $\epsilon<G(x^0)$.

\begin{remark}\label{obs:impossib-geral}
Both conditions of Proposition~\ref{prop:impossibilidade} on the steps are essential. If some block contains all the indices $1,\dots,r$, a step on that block may decrease all the objectives. If temporary increases are accepted, steps on proper blocks alone may reach a Pareto stationary point. For instance, in the example of the Introduction, starting from $x=0$, the steps $x_1=-\frac{1}{2}$ and then $x_2=-\frac{1}{2}$ yield
$$
  F(0,0)=(0,0),\qquad F(-\tfrac{1}{2},0)=(-\tfrac{3}{8},\tfrac{1}{8}),\qquad F(-\tfrac{1}{2},-\tfrac{1}{2})=(-\tfrac{1}{4},-\tfrac{1}{4}).
$$
Each step increases one of the objectives, but the final point is Pareto stationary and improves both objectives with respect to $x=0$. Thus, Proposition~\ref{prop:impossibilidade} forces steps outside the proper blocks only on monotone methods. Nonmonotone block methods are not studied in this paper.
\end{remark}

Algorithm~\ref{alg:bcd} uses $G_0(x)=G(x)$ as its only stopping criterion. The proper blocks are examined as long as one of them satisfies $G_j(x^k)>\epsilon$. When all of them satisfy $G_j(x^k)\leq\epsilon$ and still $G(x^k)>\epsilon$, the algorithm resorts to a step on block~0, that is, to a joint step in all the variables, analogous to that of Algorithm~2.1 in \cite{CDM20} with $p=1$. The computation of $G_0$ is only needed when the proper blocks are exhausted.

\begin{algorithm}[H]
\caption{Block coordinate descent with quadratic regularization for multiobjective optimization}\label{alg:bcd}
Suppose that $x^0\in\mathbb{R}^n$, $\alpha\in(0,1)$, $\theta>0$, $\epsilon\in(0,1)$, $\sigma_0\geq1$ and a partition of the variables into $q$ proper blocks are given. Set $\sigma_\ell^0:=\sigma_0$ for $\ell=0,1,\dots,q$ and $k\leftarrow0$.

\textbf{Step 1.} Set $C\leftarrow\{1,\dots,q\}$.

\textbf{Step 2.} If $C=\emptyset$, go to Step~7. Otherwise, choose $j\in C$.

\textbf{Step 3.} If $G_j(x^k)\leq\epsilon$, set $C\leftarrow C\setminus\{j\}$ and go back to Step~2.

\textbf{Step 4.} Set $\sigma\leftarrow\sigma_j^k$.

\textbf{Step 5.} If $\lambda_{\min}\bigl(B_{i,(j)}(x^k)+2\sigma I_{n_j}\bigr)<\theta$ for some $i$, set $\sigma\leftarrow2\sigma$ and repeat Step~5. Otherwise, compute an acceptable pair $(s^k_{(j)},\gamma^k)$ in the sense of Definition~\ref{def:aceitavel}, with the models \eqref{eq:modelo} defined by this $\sigma$, and set $x^{\mathrm{trial}}:=x^k+U_js^k_{(j)}$.

\textbf{Step 6.} If
$$
  f_i(x^{\mathrm{trial}})\leq f_i(x^k)-\frac{\alpha}{16\sigma}\,\epsilon^2,\qquad i=1,\dots,r,
$$
set $x^{k+1}:=x^{\mathrm{trial}}$, choose $\sigma_\ell^{k+1}\in[1,\max\{\sigma,\sigma_\ell^k\}]$ for $\ell=0,1,\dots,q$, set $k\leftarrow k+1$ and go back to Step~1. Otherwise, set $\sigma\leftarrow2\sigma$ and go back to Step~5.

\textbf{Step 7.} If $G_0(x^k)\leq\epsilon$, stop and return $x^k$ as an $\epsilon$-Pareto stationary point. Otherwise, set $j\leftarrow0$ and go to Step~4.
\end{algorithm}

\begin{remark}\label{obs:algoritmo}
By construction, $\sigma_\ell^k\geq1$ for all $\ell$ and $k$, and $\sigma$ only increases within an iteration. By Assumption~\ref{hip:matrizes},
$$
  \lambda_{\min}\big(B_{i,(j)}(x^k)+2\sigma I_{n_j}\big)\geq2\sigma-B,
$$
so that a doubling in Step~5 only occurs if $\sigma<(B+\theta)/2$. Hence, Step~5 terminates, and Lemma~\ref{lem:subproblema} guarantees the existence of the acceptable pair. When the matrices $B_{i,(j)}(x^k)$ are positive semidefinite, $\sigma\geq\theta/2$ suffices, and the regularization does not need to dominate the curvature of the model.

The rule in Step~6 includes the single nondecreasing parameter of \cite{AASS26}, given by $\sigma_\ell^{k+1}=\sigma$ for all $\ell$, and the rule with one parameter per block
$$
  \sigma_j^{k+1}=\max\{1,\sigma/\rho\},\qquad \sigma_\ell^{k+1}=\sigma_\ell^k\quad(\ell\neq j),\qquad \rho\in(1,2],
$$
which allows the parameter of each block to decrease.

Algorithm~\ref{alg:bcd} only stops at Step~7, with $G(x^k)\leq\epsilon$. As long as $G(x^k)>\epsilon$, some step is accepted at each iteration, either on a proper block with $G_j(x^k)>\epsilon$ or on block~0, after finitely many doublings of $\sigma$ (Theorem~\ref{teo:descida}).
\end{remark}

\section{Convergence analysis}\label{sec:convergencia} This section shows that each iteration of Algorithm~\ref{alg:bcd} terminates and establishes its complexity bounds.
\subsection{Sufficient decrease}
\label{sec:descida}

We now show that the decrease test in Step~6 is satisfied whenever $\sigma$ is sufficiently large and $G_j(x^k)>\epsilon$, which always holds when a step is computed, so that each iteration of Algorithm~\ref{alg:bcd} terminates after finitely many doublings of $\sigma$. Lemma~\ref{lem:decrescimo} provides an upper bound for the variation of each objective along an acceptable pair, and Lemma~\ref{lem:cota-passo} provides a lower bound for the step size when $G_j(x^k)>\epsilon$. Combined, they lead to Theorem~\ref{teo:descida}.

\begin{lemma}\label{lem:decrescimo}
Suppose that Assumptions~\ref{hip:holder} and~\ref{hip:matrizes} hold. Let $j\in\{0,1,\dots,q\}$, $\sigma>0$ and $s\in\mathbb{R}^{n_j}$ be such that $\varphi_i(s)\leq0$ for all $i=1,\dots,r$, with $\varphi_i$ defined in \eqref{eq:modelo} using $\sigma$. Then, for all $i=1,\dots,r$,
$$
  f_i(x^k+U_js)\leq f_i(x^k)-\Big(\sigma-\frac{B}{2}\Big)\|s\|^2+L_{i,j}\|s\|^{1+\beta_{i,j}} .
$$
\end{lemma}

\begin{proof}
Fix $i$ and let $g_i:=\nabla_{(j)}f_i(x^k)$ and $B_i:=B_{i,(j)}(x^k)$. By the definition of $\varphi_i$ and since $\varphi_i(s)\leq0$,
$$
  \langle g_i,s\rangle=\varphi_i(s)-\frac{1}{2}\langle B_is,s\rangle-\sigma\|s\|^2\leq-\frac{1}{2}\langle B_is,s\rangle-\sigma\|s\|^2 .
$$
By the Cauchy--Schwarz inequality and Assumption~\ref{hip:matrizes},
$$
  -\frac{1}{2}\langle B_is,s\rangle\leq\frac{1}{2}\|B_i\|\|s\|^2\leq\frac{B}{2}\|s\|^2 .
$$
Substituting into Assumption~\ref{hip:holder}, we obtain
$$
  f_i(x^k+U_js)\leq f_i(x^k)+\langle g_i,s\rangle+L_{i,j}\|s\|^{1+\beta_{i,j}}
  \leq f_i(x^k)-\Big(\sigma-\frac{B}{2}\Big)\|s\|^2+L_{i,j}\|s\|^{1+\beta_{i,j}} .
$$
\end{proof}

The proof follows that of Lemma~2 in \cite{AASS26} and involves no cross terms between the objectives, since condition \eqref{eq:aceit1} requires $\varphi_i(s)\leq0$ for each $i$, and not only for a convex combination of the models.

\begin{lemma}\label{lem:cota-passo}
Suppose that Assumption~\ref{hip:matrizes} holds. Let $j\in\{0,1,\dots,q\}$, $\sigma>0$ and $(s,\gamma)$ be an acceptable pair in the sense of Definition~\ref{def:aceitavel}, with $\varphi_i$ defined using $\sigma$. Then the following statements hold.
\begin{enumerate}
  \item[(a)] $G_j(x^k)\leq(\theta+B+2\sigma)\|s\|$.
  \item[(b)] If, in addition, $\sigma\geq(B+\theta)/2$, then $\sigma\|s\|\geq G_j(x^k)/4$.
\end{enumerate}
\end{lemma}

\begin{proof}\leavevmode
\begin{enumerate}
  \item[(a)] Let $g_i:=\nabla_{(j)}f_i(x^k)$ and $B_i:=B_{i,(j)}(x^k)$. Since $\nabla\varphi_i(s)=g_i+B_is+2\sigma s$ and $\sum_i\gamma_i=1$,
  $$
    \sum_{i=1}^r\gamma_ig_i=\sum_{i=1}^r\gamma_i\nabla\varphi_i(s)-\Big(\sum_{i=1}^r\gamma_iB_i\Big)s-2\sigma s .
  $$
  By condition \eqref{eq:aceit2}, the triangle inequality and the bound $\|\sum_i\gamma_iB_i\|\leq\sum_i\gamma_i\|B_i\|\leq B$, which follows from Assumption~\ref{hip:matrizes},
  $$
    \Big\|\sum_{i=1}^r\gamma_ig_i\Big\|\leq(\theta+B+2\sigma)\|s\| .
  $$
  Since $G_j(x^k)$ is the minimum over $\Delta^{r-1}$ and $\|\cdot\|_\infty\leq\|\cdot\|$,
  $$
    G_j(x^k)\leq\Big\|\sum_{i=1}^r\gamma_ig_i\Big\|_\infty\leq\Big\|\sum_{i=1}^r\gamma_ig_i\Big\|\leq(\theta+B+2\sigma)\|s\| .
  $$
  \item[(b)] If $\sigma\geq(B+\theta)/2$, then $\theta+B+2\sigma\leq4\sigma$, and the inequality follows from item~(a).
\end{enumerate}
\end{proof}

\begin{remark}\label{obs:lemas-descida}
Lemmas~\ref{lem:decrescimo} and~\ref{lem:cota-passo} extend Lemmas~2 and~3 of \cite{AASS26}. In Lemma~\ref{lem:cota-passo}, the inequality in item~(a) holds for the $\gamma$ of any acceptable pair, because $G_j$ is a minimum over $\Delta^{r-1}$, and it involves no constants depending on $n_j$ or $r$.
\end{remark}

To quantify when the test in Step~6 is satisfied, we define, for $i=1,\dots,r$, $j\in\{0,1,\dots,q\}$ and $t>0$,
\begin{equation}\label{eq:sigmaij}
\sigma_{i,j}(t):=\bigg[\Big(\frac{B}{2(1-\alpha)}+\frac{L_{i,j}}{1-\alpha}\Big)\Big(\frac{t}{4}\Big)^{\beta_{i,j}-1}\bigg]^{1/\beta_{i,j}},
\qquad \sigma_j(t):=\max_{1 \leq i \leq r}\sigma_{i,j}(t).
\end{equation}
Since $\beta_{i,j}\in(0,1]$, we have $\sigma_{i,j}(\epsilon)=\epsilon^{(\beta_{i,j}-1)/\beta_{i,j}}\sigma_{i,j}(1)\geq\sigma_{i,j}(1)$ for all $\epsilon\in(0,1)$.

\begin{theorem}\label{teo:descida}
Suppose that Assumptions~\ref{hip:holder} and~\ref{hip:matrizes} hold. Let $j\in\{0,1,\dots,q\}$, $\epsilon\in(0,1)$ and $\sigma\geq\max\{1,(B+\theta)/2,\sigma_j(\epsilon)\}$, and let $(s,\gamma)$ be an acceptable pair in the sense of Definition~\ref{def:aceitavel}, with $\varphi_i$ defined using $\sigma$. If $G_j(x^k)>\epsilon$, then, for all $i=1,\dots,r$,
\begin{equation}\label{eq:descida}
  f_i(x^k+U_js)\leq f_i(x^k)-\alpha\sigma\|s\|^2
  \qquad\text{and}\qquad
  f_i(x^k+U_js)\leq f_i(x^k)-\frac{\alpha}{16\sigma}\epsilon^2 .
\end{equation}
In particular, $x^{\mathrm{trial}}=x^k+U_js$ satisfies the test in Step~6 of Algorithm~\ref{alg:bcd}.
\end{theorem}

The proof, which adapts that of Theorem~4 in \cite{AASS26}, is given in Appendix~\ref{ap:provas}. The multiobjective case enters at two points: in the bound $\sigma\|s\|>\epsilon/4$, obtained from Lemma~\ref{lem:cota-passo}(b) through $\sum_i\gamma_i\nabla_{(j)}f_i(x^k)$, and in the maximum over $i$ that defines $\sigma_j(t)$.

\subsection{Complexity bounds}
\label{sec:complexidade}

In this subsection, we establish complexity bounds for Algorithm~\ref{alg:bcd}. We call \emph{outer iteration} each increment of $k$, that is, each step accepted in Step~6. Rejected trials are counted only in the number of evaluations of $F$. Lemma~\ref{lem:decr-uniforme} guarantees a uniform decrease per accepted step, from which follow the iteration bounds (Theorem~\ref{teo:iteracoes}), the evaluation bounds (Theorem~\ref{teo:avaliacoes}) and convergence with decreasing tolerances (Corollary~\ref{cor:convergencia}). Theorem~\ref{teo:melhor} deals with alternative acceptance tests, and Remark~\ref{obs:casos} deals with the case $r=1$ and with the exclusive use of block~0.

\begin{lemma}\label{lem:decr-uniforme}
Suppose that Assumptions~\ref{hip:holder} and~\ref{hip:matrizes} hold, and let
$$
  \beta:=\min_{i,j}\beta_{i,j},\qquad c:=\max\Bigl\{\sigma_0,\ B+\theta,\ 2\max_{i,j}\sigma_{i,j}(1)\Bigr\},\qquad \sigma_{\max}:=c\,\epsilon^{(\beta-1)/\beta},
$$
where the maxima and minima are taken over $i\in\{1,\dots,r\}$ and $j\in\{0,1,\dots,q\}$, and let $\eta:=16c$. Then, as long as Algorithm~\ref{alg:bcd} does not stop, the following statements hold.
\begin{enumerate}
  \item[(a)] $\sigma_\ell^k\leq\sigma_{\max}$ for all $\ell\in\{0,1,\dots,q\}$ and $k$, and every value of $\sigma$ used in Step~5 satisfies $\sigma\leq\sigma_{\max}$.
  \item[(b)] For all $k$ and all $i=1,\dots,r$,
  $$
    f_i(x^{k+1})\leq f_i(x^k)-\frac{\alpha}{\eta}\,\epsilon^{(\beta+1)/\beta}.
  $$
\end{enumerate}
\end{lemma}

The proof, which follows those of Lemma~5 and Theorem~8 in \cite{AASS26}, is given in Appendix~\ref{ap:provas}. Compared with \cite{AASS26}, it also handles the doublings in Step~5 and the general update rule in Step~6, which includes one parameter per block.

\begin{theorem}\label{teo:iteracoes}
Suppose that Assumptions~\ref{hip:holder} to~\ref{hip:limitacao} hold, and let $\beta$ and $\eta$ be as in Lemma~\ref{lem:decr-uniforme}. Then Algorithm~\ref{alg:bcd} stops at a point $x^K$ with $G(x^K)\leq\epsilon$, where the number $K$ of outer iterations satisfies
$$
  K\leq\frac{\eta\,\big(f_{i_0}(x^0)-f_{\mathrm{low}}\big)}{\alpha}\,\epsilon^{-(\beta+1)/\beta}.
$$
\end{theorem}

The proof, which follows that of Theorem~8 in \cite{AASS26} applied to the objective $f_{i_0}$ of Assumption~\ref{hip:limitacao}, is given in Appendix~\ref{ap:provas}. If more than one objective satisfies Assumption~\ref{hip:limitacao}, the bound holds for each of them, and one may take the smallest of the right-hand sides.

The order $O(\epsilon^{-(\beta+1)/\beta})$ in Theorem~\ref{teo:iteracoes} is the one obtained in \cite{CGT17} for quadratic regularization and in \cite{CDM20} with $p=1$. In the Lipschitz case, $\beta=1$, it is $O(\epsilon^{-2})$ and, for the variant with the Euclidean norm in \eqref{eq:Gj}, it is optimal. Indeed, with $f_1=\dots=f_r$, the measure $G$ reduces to $\|\nabla f_1(x)\|$. The hard functions of \cite{CDHS20} have Lipschitz continuous gradients and are bounded below, and therefore satisfy Assumptions~\ref{hip:holder} to~\ref{hip:limitacao} with $\beta=1$ and $B_{i,(j)}\equiv0$. For them, the number of queries to an oracle that returns the function value and the gradient at a point is at least of the order of $\epsilon^{-2}$, both for deterministic algorithms \cite[Theorem~1]{CDHS20} and for randomized algorithms \cite[Theorem~2]{CDHS20}. For deterministic algorithms, this bound was extended in \cite[Theorem~4.12]{Sam26} to multiobjective problems with distinct objectives. The dimension of these hard functions grows with $\epsilon^{-2}$, so that optimality refers to the worst case over all dimensions. For $\beta=1$, we have $\sigma_{\max}=c$ and, by the count \eqref{eq:contagem}, established below in the proof of Theorem~\ref{teo:avaliacoes}, and by Lemma~\ref{lem:decr-uniforme}(a), each iteration of Algorithm~\ref{alg:bcd} uses at most $2+\log_2c$ queries, one at $x^k$ and the others at the trial points, a number that does not depend on $\epsilon$, $n$ or $q$. Hence, the lower bound also holds for the number of iterations. With the maximum norm, $G(x)\leq\epsilon$ only guarantees that the Euclidean measure is at most $\sqrt n\,\epsilon$ (Lemma~\ref{lem:prop-Gj}(d)). For $\beta<1$, we are not aware of a lower bound in the multiobjective case.

\begin{corollary}\label{cor:convergencia}
Suppose that Assumptions~\ref{hip:holder} to~\ref{hip:limitacao} hold. Let $\{\epsilon_k\}\subset(0,1)$ be a sequence with $\epsilon_k\to0$. Set $\bar x^0:=x^0$ and, for $k\geq0$, let $\bar x^{k+1}$ be the point returned by Algorithm~\ref{alg:bcd} with tolerance $\epsilon_k$ and initial point $\bar x^k$. Then the following statements hold.
\begin{enumerate}
  \item[(a)] Each run terminates after finitely many outer iterations, and $G(\bar x^{k+1})\leq\epsilon_k$ for all $k$.
  \item[(b)] Every accumulation point of $\{\bar x^k\}$ is Pareto stationary.
  \item[(c)] If the level set $\mathcal{L}_0$ of Assumption~\ref{hip:limitacao} is bounded, then $\{\bar x^k\}$ has at least one accumulation point.
\end{enumerate}
\end{corollary}

\begin{proof}\leavevmode
\begin{enumerate}
  \item[(a)] Since each accepted step decreases $f_{i_0}$, we have $f_{i_0}(\bar x^k)\leq f_{i_0}(x^0)$ and, therefore,
  $$
    \{x\in\mathbb{R}^n:\ f_{i_0}(x)\leq f_{i_0}(\bar x^k)\}\subset\mathcal{L}_0 .
  $$
  Thus, Assumption~\ref{hip:limitacao} holds with the same $f_{\mathrm{low}}$ for the initial point $\bar x^k$, and the result follows from Theorem~\ref{teo:iteracoes}.
  \item[(b)] Let $\bar x^*$ be an accumulation point of $\{\bar x^k\}$ and $\{\bar x^{k_l}\}$, with $k_l\geq1$, a subsequence converging to $\bar x^*$. By item~(a), $G(\bar x^{k_l})\leq\epsilon_{k_l-1}\to0$. By the continuity of $G$ (Lemma~\ref{lem:prop-Gj}(c)),
  $$
    G(\bar x^*)=\lim_{l\to\infty}G(\bar x^{k_l})=0,
  $$
  and, therefore, $\bar x^*$ is Pareto stationary.
  \item[(c)] The set $\mathcal{L}_0$ is closed, by the continuity of $f_{i_0}$, and contains $\{\bar x^k\}$, by the proof of item~(a). If $\mathcal{L}_0$ is bounded, then it is compact and, by the Bolzano--Weierstrass theorem, $\{\bar x^k\}$ has a convergent subsequence, whose limit is an accumulation point.
\end{enumerate}
\end{proof}

If $\{\epsilon_k\}$ is nonincreasing, Corollary~\ref{cor:convergencia} can be viewed as a single run of a variant of Algorithm~\ref{alg:bcd} in which, whenever Step~7 verifies $G(x^k)\leq\epsilon_l$ for the current tolerance $\epsilon_l$, the tolerance is changed to $\epsilon_{l+1}$ and the algorithm continues from the same point. The parameters $\sigma_\ell$ may either be reset to $\sigma_0$ or kept. In the latter case, the proof of Lemma~\ref{lem:decr-uniforme}(a) remains valid, since $\sigma_{\max}$ does not decrease when $\epsilon$ decreases. For this variant, $\liminf_kG(x^k)=0$, and every accumulation point of the subsequence of points at which the tolerance is updated is Pareto stationary. Whether $\lim_kG(x^k)=0$ remains open, since, by Proposition~\ref{prop:impossibilidade}, the measures $G_j$ may be small while $G$ is not.

To measure the cost in gradients, we count the computation of $\nabla f_1(x),\dots,\nabla f_r(x)$ at a point as one \emph{gradient equivalent}, and that of $\nabla_{(j)}f_1(x),\dots,\nabla_{(j)}f_r(x)$ as $n_j/n$ gradient equivalents, without including the cost of the matrices $B_{i,(j)}$.

\begin{theorem}\label{teo:avaliacoes}
Suppose that Assumptions~\ref{hip:holder} to~\ref{hip:limitacao} hold, let $\sigma_{\max}$ and $\eta$ be as in Lemma~\ref{lem:decr-uniforme}, and let
$$
  K_\epsilon:=\frac{\eta\,\bigl(f_{i_0}(x^0)-f_{\mathrm{low}}\bigr)}{\alpha}\,\epsilon^{-(\beta+1)/\beta}.
$$
Then, until it stops, Algorithm~\ref{alg:bcd} satisfies the following bounds.
\begin{enumerate}
  \item[(a)] The number of evaluations of $F$ is at most $1+K_\epsilon\bigl(1+\log_2\sigma_{\max}\bigr)$.
  \item[(b)] If, at each step accepted on block $j$ with parameter $\sigma$, the rule in Step~6 satisfies $\sigma_j^{k+1}\geq\sigma/2$ and $\sigma_\ell^{k+1}\geq\sigma_\ell^k$ for $\ell\neq j$, then the number of evaluations of $F$ is at most $1+2K_\epsilon+(q+1)\log_2(\sigma_{\max}/\sigma_0)$.
  \item[(c)] If the rule in Step~6 assigns the same value $\sigma^{k+1}$ to all blocks, that is, $\sigma_\ell^{k+1}=\sigma^{k+1}$ for all $\ell$, then the number of evaluations of $F$ is at most $1+2K_\epsilon+\log_2(\sigma_{\max}/\sigma_0)$ when $\sigma^{k+1}\geq\sigma/2$, and at most $1+K_\epsilon+\log_2(\sigma_{\max}/\sigma_0)$ when $\sigma^{k+1}\geq\sigma$.
  \item[(d)] The measures $G_j$, with $j\geq1$, are computed at most $q\,(K_\epsilon+1)$ times, and $G_0$ at most $K_\epsilon+1$ times.
  \item[(e)] The number of gradient equivalents is at most $K_\epsilon+1$.
\end{enumerate}
\end{theorem}

\begin{proof}
Let $K\leq K_\epsilon$ be the number of outer iterations (Theorem~\ref{teo:iteracoes}) and, for $k=0,\dots,K-1$, let $m_k$ be the number of doublings of $\sigma$, in Steps~5 and~6, during iteration $k$. Besides $F(x^0)$, each execution of Step~6 evaluates $F$ once, at $x^{\mathrm{trial}}$. The doublings in Step~5 do not evaluate $F$, and the last iteration ends at Step~7 without evaluating $F$. Hence, the total number of evaluations is at most
\begin{equation}\label{eq:contagem}
  1+K+\sum_{k=0}^{K-1}m_k.
\end{equation}
Items~(a) and~(c), which adapt the proof of Theorem~7 in \cite{AASS26}, are proved in Appendix~\ref{ap:provas}. We prove items~(b), (d) and~(e) here.
\begin{enumerate}
    \item[(b)] Let $\Phi_k:=\sum_{\ell=0}^q\log_2\sigma_\ell^k$. At iteration $k$, with block $j$, the accepted parameter is $\sigma=2^{m_k}\sigma_j^k$. By the assumption of item~(b),
  $$
    \log_2\sigma_j^{k+1}\geq\log_2\sigma-1=\log_2\sigma_j^k+m_k-1
    \qquad\text{and}\qquad
    \log_2\sigma_\ell^{k+1}\geq\log_2\sigma_\ell^k,\quad \ell\neq j .
  $$
  Hence, $\Phi_{k+1}-\Phi_k\geq m_k-1$. Summing over $k=0,\dots,K-1$ and using $\sigma_\ell^0=\sigma_0$ and $\sigma_\ell^K\leq\sigma_{\max}$ (Lemma~\ref{lem:decr-uniforme}(a)), we obtain
  $$
    \sum_{k=0}^{K-1}m_k\leq K+\Phi_K-\Phi_0\leq K+(q+1)\log_2(\sigma_{\max}/\sigma_0),
  $$
  and item~(b) follows from \eqref{eq:contagem} and $K\leq K_\epsilon$.
    \item[(d)] At each iteration $k=0,1,\dots,K$, including the last one, each proper block is tested in Step~3 at most once, since it is either removed from $C$ or the step is accepted and the iteration ends. Step~7 is executed at most once per iteration. Since there are $K+1\leq K_\epsilon+1$ iterations, item~(d) follows.
    \item[(e)] Fix an iteration $k\in\{0,1,\dots,K\}$. In Step~3, the test of block $j\in\{1,\dots,q\}$ requires $\nabla_{(j)}f_1(x^k),\dots,\nabla_{(j)}f_r(x^k)$, which cost $n_j/n$ gradient equivalents and, by item~(d), each block is tested at most once. Step~7 is reached only when $C=\emptyset$, that is, after all proper blocks have been tested in Step~3 at $x^k$. At that point, the partial gradients $\nabla_{(j)}f_i(x^k)$, $j=1,\dots,q$, are available and form $\nabla f_i(x^k)$, so the computation of $G_0(x^k)$ requires no new gradients. Steps~5 and~6, including the rejected trials, do not require new gradients either, since they only change $\sigma$ and use the gradients already computed at $x^k$. Since $n_1+\dots+n_q=n$, the cost of iteration $k$ is at most
  $$
    \sum_{j=1}^q\frac{n_j}{n}=1 .
  $$
  Summing over the $K+1\leq K_\epsilon+1$ iterations, we obtain the bound $K_\epsilon+1$.
\end{enumerate}
\end{proof}

Both rules of Remark~\ref{obs:algoritmo} are covered by Theorem~\ref{teo:avaliacoes}. The rule with one parameter per block satisfies the assumption of item~(b), since $\max\{1,\sigma/\rho\}\geq\sigma/2$ for $\rho\leq2$, and the rule with a single nondecreasing parameter satisfies the second part of item~(c), with a bound of the same form as that of Theorem~7 in \cite{AASS26}. In these cases, the number of evaluations of $F$ is $O(\epsilon^{-(\beta+1)/\beta})$, whereas the general bound in item~(a) has an additional factor $\log_2\sigma_{\max}=O(\log(1/\epsilon))$.

Since a method that updates all the variables uses at least one gradient equivalent per iteration, item~(e) shows that, in the worst case, Algorithm~\ref{alg:bcd} uses at most as many gradients as such a method with the same number of iterations. Any advantage of block descent must therefore come from the cost observed in practice, not from this bound; Remark~\ref{obs:porbloco} gives bounds that depend on the constants of each block.

\begin{remark}\label{obs:porbloco}
If $\sigma_j^{k+1}\leq\sigma$ and $\sigma_\ell^{k+1}=\sigma_\ell^k$ for $\ell\neq j$, as in the rule with one parameter per block, the bound of Lemma~\ref{lem:decr-uniforme}(a) holds blockwise. Let
$$
  \bar\sigma_\ell(\epsilon):=\max\{(B+\theta)/2,\,\sigma_\ell(\epsilon)\},\qquad \hat\sigma_\ell:=\max\{\sigma_0,\,2\bar\sigma_\ell(\epsilon)\}.
$$
By \eqref{eq:cota-sigmaj}, obtained in the proof of Lemma~\ref{lem:decr-uniforme} in Appendix~\ref{ap:provas}, we have $\hat\sigma_\ell\leq\sigma_{\max}$, and every $\sigma$ used with block $\ell$ satisfies $\sigma\leq\hat\sigma_\ell$. Indeed, $\sigma_\ell^k$ only changes when block $\ell$ is chosen, and it is then doubled only if $\sigma<\bar\sigma_\ell(\epsilon)$. It suffices to apply to each block the induction in the proof of Lemma~\ref{lem:decr-uniforme}(a). As in the proof of Lemma~\ref{lem:decr-uniforme}(b), a step accepted on block $j$ satisfies
$$
  f_i(x^{k+1})\leq f_i(x^k)-\frac{\alpha\epsilon^2}{16\hat\sigma_j},\qquad i=1,\dots,r,
$$
with at most $\log_2\hat\sigma_j$ doublings. Thus, the decrease and the doublings in an iteration on block $j$ depend, among the smoothness constants, only on the constants $L_{i,j}$ and $\beta_{i,j}$ of that block. Moreover, let $N_\ell$ be the number of steps accepted on block $\ell$ until the algorithm stops. Applying the inequality above to $f_{i_0}$ and summing over all iterations, as in the proof of Theorem~\ref{teo:iteracoes}, we obtain
$$
  \sum_{\ell=0}^q\frac{N_\ell}{\hat\sigma_\ell}\leq\frac{16\,\bigl(f_{i_0}(x^0)-f_{\mathrm{low}}\bigr)}{\alpha\,\epsilon^2}.
$$
In particular, $N_\ell\leq16\,\hat\sigma_\ell\bigl(f_{i_0}(x^0)-f_{\mathrm{low}}\bigr)/(\alpha\epsilon^2)$ for each $\ell$, so that a block with small constants can be chosen only a few times. Since $\hat\sigma_\ell\leq\sigma_{\max}$ and $K=\sum_\ell N_\ell$, the bound of Theorem~\ref{teo:iteracoes} follows from this weighted bound. The weighted bound reduces to the uniform one when a block with $\hat\sigma_j$ close to $\sigma_{\max}$ is chosen at every iteration, but it shows that the number of steps on each block is controlled by the constants of that block. In contrast, with a single nondecreasing parameter, a large $\sigma$ obtained on one block is used on all subsequent ones. This suggests an advantage of the per-block rule when the constants vary across blocks.
\end{remark}

The following result deals with two acceptance tests alternative to the one in Step~6.

\begin{theorem}\label{teo:melhor}
Suppose that Assumptions~\ref{hip:holder} to~\ref{hip:limitacao} hold, and let $\beta$, $c$ and $\sigma_{\max}$ be as in Lemma~\ref{lem:decr-uniforme}. Then the following statements hold.
\begin{enumerate}
  \item[(a)] If, in Step~6, the acceptance test is replaced by the test of \cite{AASS26}
  \begin{equation}\label{eq:teste-misto}
    f_i(x^{\mathrm{trial}})\leq f_i(x^k)-\min\Big\{\frac{\alpha}{16\sigma}\,\epsilon^2,\ \alpha\|s^k_{(j)}\|^2\Big\},\qquad i=1,\dots,r,
  \end{equation}
  then Algorithm~\ref{alg:bcd} stops at a point $x^K$ with $G(x^K)\leq\epsilon$, where the number $K$ of outer iterations satisfies $K\leq\hat K_\epsilon$, with
  $$
    \hat K_\epsilon:=\frac{\hat\eta\,\big(f_{i_0}(x^0)-f_{\mathrm{low}}\big)}{\alpha}\,\epsilon^{-2/\beta}
    \qquad\text{and}\qquad
    \hat\eta:=16c^2 .
  $$
  \item[(b)] If the test is replaced by
  \begin{equation}\label{eq:teste-sigma}
    f_i(x^{\mathrm{trial}})\leq f_i(x^k)-\alpha\sigma\|s^k_{(j)}\|^2,\qquad i=1,\dots,r,
  \end{equation}
  and if $\sigma_0\geq(B+\theta)/2$ and $\sigma_\ell^{k+1}\geq(B+\theta)/2$ in the rule of Step~6, then the bound of Theorem~\ref{teo:iteracoes} holds.
\end{enumerate}
In both cases, the bounds on the number of evaluations of $F$ in items~(a)--(c) of Theorem~\ref{teo:avaliacoes} hold with $K_\epsilon$ replaced by the corresponding iteration bound.
\end{theorem}

The proof, which follows those of Theorems~6 and~7 and of Remark~2 in \cite{AASS26}, is given in Appendix~\ref{ap:provas}.

Although test \eqref{eq:teste-misto} of \cite{AASS26} is the least demanding of the three, its bound has a strictly worse order when $\beta<1$, since $(\beta+1)/\beta<2/\beta$. The orders coincide when $\beta=1$. The test in Step~6 is that of Theorem~8 in \cite{AASS26} and, up to the constant, that of \cite{CDM20}. Test \eqref{eq:teste-sigma}, proposed in \cite[Remark~2]{AASS26}, does not depend on $\epsilon$, but requires $\sigma\geq(B+\theta)/2$ throughout the algorithm.

\begin{remark}\label{obs:casos}
For $r=1$, we have $G_j(x)=\|\nabla_{(j)}f_1(x)\|_\infty$, conditions \eqref{eq:aceit1} and \eqref{eq:aceit2} are conditions (2) and (3) of \cite{AASS26}, and, by Lemma~\ref{lem:prop-Gj}(b), Algorithm~\ref{alg:bcd} stops at Step~7 without using block~0. Suppose, in addition, that $\sigma_\ell^{k+1}=\sigma$ for all $\ell$, that the test is \eqref{eq:teste-misto}, that the matrices $B_{(j)}(x^k)$ are positive semidefinite, and that $\theta\leq2$, so that Step~5 never doubles $\sigma$. Then Algorithm~\ref{alg:bcd} generates the same iterates as Algorithm~1 of \cite{AASS26} with $f_{\mathrm{est}}=-\infty$. In this case, Assumption~\ref{hip:holder} is only needed for $j\geq1$, and the bounds of Theorems~\ref{teo:iteracoes}, \ref{teo:avaliacoes} and~\ref{teo:melhor} have the same order as those of Theorems~6 to~8 of \cite{AASS26}, with different constants. With the per-block rule, the iterates differ, but the orders are the same.

In the other extreme case, if Step~1 sets $C\leftarrow\emptyset$, only block~0 is used. The proofs hold without change, and each iteration uses a single gradient equivalent. With $B_{i,(0)}\equiv0$, subproblem \eqref{eq:subproblema} is that of \cite{CDM20} without constraints, with $p=1$ and a linear model, and the order in Theorem~\ref{teo:iteracoes} coincides with that of Theorem~2.3 in \cite{CDM20}, since $(p+\beta)/(p+\beta-1)=(1+\beta)/\beta$ for $p=1$. However, the iterates differ, since \cite{CDM20} requires the decrease $\alpha\epsilon^2/(36\sigma)$ and checks the stopping criterion at the trial point, in the Euclidean norm. When $B_{i,(0)}\not\equiv0$ and $\beta<1$, the step on block~0 is an extension of the method of \cite{CDM20}.
\end{remark}

\section{Conclusions}\label{sec:conclusoes}

We proposed a block coordinate descent method for unconstrained, possibly nonconvex, multiobjective optimization with H\"older continuous gradients, which extends the method of \cite{AASS26} and admits quadratic models on each block, including second-order information. We showed that this extension is not straightforward. With two or more objectives, the blockwise stationarity measures do not guarantee Pareto stationarity, and there are problems with Lipschitz continuous gradients for which no method that uses only proper blocks and never increases any objective reaches $\epsilon$-Pareto stationarity for small $\epsilon$ (Proposition~\ref{prop:impossibilidade}). For this reason, the method resorts to a step in all the variables only when the measures of all proper blocks do not exceed the tolerance and the point is not yet $\epsilon$-Pareto stationary, and the stopping criterion guarantees that the final point is indeed $\epsilon$-Pareto stationary.

We proved a bound of $O(\epsilon^{-(\beta+1)/\beta})$ iterations which, in the Lipschitz case and for the variant with the Euclidean norm, is optimal in the worst case over all dimensions (Section~\ref{sec:complexidade}). We further proved bounds on the numbers of function and gradient evaluations, valid with one regularization parameter per block, a bound on the number of steps on each block that depends on the constants of that block, and the Pareto stationarity of the accumulation points with decreasing tolerances. With a single objective, we recover the orders of \cite{AASS26}. Using only the step in all the variables, we recover the order of \cite{CDM20} with $p=1$.

As future work, we highlight implementations of the second-order models that exploit sparsity, higher-order regularizations on each block, as in \cite{CDM20}, block steps with curvature information of Barzilai--Borwein type \cite{morovati2016,chen2023}, universal variants that do not require knowledge of the H\"older constants, as in \cite{PG25}, and methods with nonmonotone acceptance, which may dispense with the step in all the variables in some cases (Remark~\ref{obs:impossib-geral}). A numerical study of the method will be reported in a subsequent version of this work.

\backmatter

\begin{appendices}

 \section{Proofs adapted from the scalar case}\label{ap:provas}

This appendix collects the proofs that adapt those of \cite{AASS26} to the multiobjective case. The corresponding result of \cite{AASS26} and the points at which the multiobjective case enters are indicated in the main text, after each statement.

\begin{proof}[Proof of Theorem~\ref{teo:descida}]
Fix $i$. By Lemma~\ref{lem:cota-passo}(b), which applies because $\sigma\geq(B+\theta)/2$, and since $G_j(x^k)>\epsilon$,
\begin{equation}\label{eq:passo-min}
  \sigma\|s\|\geq\frac{G_j(x^k)}{4}>\frac{\epsilon}{4}.
\end{equation}
In particular, $s\neq0$. Since $(s,\gamma)$ is acceptable, Lemma~\ref{lem:decrescimo} applies, and the first inequality in \eqref{eq:descida} holds if
$$
  \Big(\sigma-\frac{B}{2}\Big)\|s\|^2-L_{i,j}\|s\|^{1+\beta_{i,j}}\geq\alpha\sigma\|s\|^2,
  \quad\text{that is, if}\quad
  \sigma\geq\frac{B}{2(1-\alpha)}+\frac{L_{i,j}}{1-\alpha}\|s\|^{\beta_{i,j}-1}.
$$
Since $\beta_{i,j}-1\leq0$, \eqref{eq:passo-min} gives $\|s\|^{\beta_{i,j}-1}\leq(\epsilon/4)^{\beta_{i,j}-1}\sigma^{1-\beta_{i,j}}$. Moreover, since $\epsilon<1$ and $\sigma\geq1$, we have $(\epsilon/4)^{\beta_{i,j}-1}\sigma^{1-\beta_{i,j}}\geq1$. Therefore,
\begin{align*}
  \frac{B}{2(1-\alpha)}+\frac{L_{i,j}}{1-\alpha}\|s\|^{\beta_{i,j}-1}
  &\leq\Big(\frac{B}{2(1-\alpha)}+\frac{L_{i,j}}{1-\alpha}\Big)\Big(\frac{\epsilon}{4}\Big)^{\beta_{i,j}-1}\sigma^{1-\beta_{i,j}}\\
  &=\sigma_{i,j}(\epsilon)^{\beta_{i,j}}\,\sigma^{1-\beta_{i,j}}
  \leq\sigma,
\end{align*}
where the equality follows from \eqref{eq:sigmaij} and the last inequality uses $\sigma\geq\sigma_j(\epsilon)\geq\sigma_{i,j}(\epsilon)$. This proves the first inequality in \eqref{eq:descida}. The second one follows from it and from \eqref{eq:passo-min}, since
$$
  \alpha\sigma\|s\|^2=\frac{\alpha(\sigma\|s\|)^2}{\sigma}\geq\frac{\alpha}{16\sigma}\,\epsilon^2.
$$
Since $i$ is arbitrary, \eqref{eq:descida} holds for all $i=1,\dots,r$, and the second inequality in \eqref{eq:descida} is the test in Step~6.
\end{proof}

\begin{proof}[Proof of Lemma~\ref{lem:decr-uniforme}]
Let $j$ be the block for which Step~4 is executed at iteration $k$. This block is unique, since the loop of Steps~5 and~6 only ends when a step is accepted. By Steps~3 and~7, we have $G_j(x^k)>\epsilon$. Define $\bar\sigma_j(\epsilon):=\max\{(B+\theta)/2,\sigma_j(\epsilon)\}$. Since $0<\epsilon<1$, $0<\beta\leq\beta_{i,j}\leq1$ and $t\mapsto(t-1)/t$ is increasing on $(0,1]$, we have $\epsilon^{(\beta_{i,j}-1)/\beta_{i,j}}\leq\epsilon^{(\beta-1)/\beta}$. Since $\sigma_{i,j}(\epsilon)=\epsilon^{(\beta_{i,j}-1)/\beta_{i,j}}\sigma_{i,j}(1)$ and $\epsilon^{(\beta-1)/\beta}\geq1$, it follows from the definition of $c$ that
\begin{equation}\label{eq:cota-sigmaj}
\bar\sigma_j(\epsilon)\leq\max\Bigl\{\frac{B+\theta}{2},\ \epsilon^{(\beta-1)/\beta}\max_{i,j}\sigma_{i,j}(1)\Bigr\}\leq\frac{\sigma_{\max}}{2}.
\end{equation}
Moreover, $\sigma_{\max}\geq c\geq\max\{\sigma_0,B+\theta\}$ and $\sigma_0\geq1$.
\begin{enumerate}
\item[(a)] By induction on $k$. For $k=0$, $\sigma_\ell^0=\sigma_0\leq\sigma_{\max}$ for all $\ell$. Suppose that $\sigma_\ell^k\leq\sigma_{\max}$ for all $\ell$. At iteration $k$, the first value used is $\sigma=\sigma_j^k\leq\sigma_{\max}$. A doubling in Step~5 only occurs if $\sigma<(B+\theta)/2$ (Remark~\ref{obs:algoritmo}). A rejection in Step~6 only occurs if $\sigma<\bar\sigma_j(\epsilon)$, since, by Theorem~\ref{teo:descida}, the test in Step~6 is satisfied whenever $\sigma\geq\bar\sigma_j(\epsilon)$, because $\sigma\geq\sigma_j^k\geq1$. In both cases, $\sigma<\bar\sigma_j(\epsilon)$, and the next value, $2\sigma$, is smaller than $2\bar\sigma_j(\epsilon)\leq\sigma_{\max}$, by \eqref{eq:cota-sigmaj}. Therefore, every $\sigma$ used at iteration $k$, including the accepted one, satisfies $\sigma\leq\sigma_{\max}$. By the rule in Step~6, $\sigma_\ell^{k+1}\leq\max\{\sigma,\sigma_\ell^k\}\leq\sigma_{\max}$ for all $\ell$.

\item[(b)] Let $\sigma$ be the parameter of the accepted step, with $\sigma\leq\sigma_{\max}$ by item~(a). By the test in Step~6,
\begin{align*}
f_i(x^{k+1})&\leq f_i(x^k)-\frac{\alpha\epsilon^2}{16\sigma}\leq f_i(x^k)-\frac{\alpha\epsilon^2}{16\sigma_{\max}}\\
&= f_i(x^k)-\frac{\alpha}{16c}\,\epsilon^{2-(\beta-1)/\beta}= f_i(x^k)-\frac{\alpha}{\eta}\,\epsilon^{(\beta+1)/\beta}.
\end{align*}
\end{enumerate}
\end{proof}

\begin{proof}[Proof of Theorem~\ref{teo:iteracoes}]
As long as Algorithm~\ref{alg:bcd} does not stop, some step is accepted at each iteration $k$ (Remark~\ref{obs:algoritmo}), and Lemma~\ref{lem:decr-uniforme}(b) gives
$$
  f_i(x^{k+1})\leq f_i(x^k)-\frac{\alpha}{\eta}\,\epsilon^{(\beta+1)/\beta},\qquad i=1,\dots,r.
$$
In particular, each accepted step decreases $f_{i_0}$, so that $x^k\in\mathcal{L}_0$ for all $k$, and Assumption~\ref{hip:limitacao} gives $f_{i_0}(x^k)\geq f_{\mathrm{low}}$. Suppose that the algorithm does not stop at iterations $0,1,\dots,K-1$. Summing the inequalities above for $i=i_0$,
$$
  f_{\mathrm{low}}\leq f_{i_0}(x^K)\leq f_{i_0}(x^0)-K\,\frac{\alpha}{\eta}\,\epsilon^{(\beta+1)/\beta},
$$
from which the bound on $K$ follows. Hence, the number of outer iterations is finite. Since the algorithm only stops at Step~7, with $G(x^K)\leq\epsilon$ (Remark~\ref{obs:algoritmo}), the theorem is proved.
\end{proof}

\begin{proof}[Proof of items~(a) and~(c) of Theorem~\ref{teo:avaliacoes}]
We use the notation $m_k$ and the count \eqref{eq:contagem} from the proof of the other items.
\begin{enumerate}
\item[(a)] At iteration $k$, the first value used is $\sigma_j^k\geq1$, and after $m_k$ doublings the value used is $2^{m_k}\sigma_j^k\leq\sigma_{\max}$, by Lemma~\ref{lem:decr-uniforme}(a). Hence, $m_k\leq\log_2\sigma_{\max}$, and item~(a) follows from \eqref{eq:contagem} and $K\leq K_\epsilon$.

\item[(c)] Since $\sigma_\ell^0=\sigma_0$ for all $\ell$, the rule in item~(c) gives $\sigma_\ell^k=\sigma^k$ for all $\ell$ and $k$, with $\sigma^0:=\sigma_0$. At iteration $k$, the accepted parameter is $\sigma=2^{m_k}\sigma^k$. If $\sigma^{k+1}\geq\sigma/2$, then $\log_2\sigma^{k+1}\geq\log_2\sigma^k+m_k-1$ and, summing as in the proof of item~(b),
$$
\sum_{k=0}^{K-1}m_k\leq K+\log_2(\sigma^K/\sigma_0)\leq K+\log_2(\sigma_{\max}/\sigma_0).
$$
If $\sigma^{k+1}\geq\sigma$, then $\log_2\sigma^{k+1}\geq\log_2\sigma^k+m_k$ and, in the same way, $\sum_{k=0}^{K-1}m_k\leq\log_2(\sigma_{\max}/\sigma_0)$. In both cases, item~(c) follows from \eqref{eq:contagem} and $K\leq K_\epsilon$.
\end{enumerate}
\end{proof}

\begin{proof}[Proof of Theorem~\ref{teo:melhor}]
By Theorem~\ref{teo:descida}, both tests are satisfied whenever $\sigma\geq\bar\sigma_j(\epsilon)$. Test \eqref{eq:teste-misto} is satisfied because the decrease it requires does not exceed $\frac{\alpha}{16\sigma}\epsilon^2$, and test \eqref{eq:teste-sigma} is satisfied by the first inequality in \eqref{eq:descida}. Hence, the proof of Lemma~\ref{lem:decr-uniforme}(a) holds without change, and every $\sigma$ used satisfies $\sigma\leq\sigma_{\max}$.
\begin{enumerate}
\item[(a)] Let $s:=s^k_{(j)}$ be the accepted step, with parameter $\sigma\leq\sigma_{\max}$. By Lemma~\ref{lem:cota-passo}(a) and since $B+\theta\leq\sigma_{\max}$,
$$
\epsilon<G_j(x^k)\leq(\theta+B+2\sigma)\|s\|\leq3\sigma_{\max}\|s\|<4\sigma_{\max}\|s\|.
$$
Since $\sigma\leq\sigma_{\max}$ and $\sigma_{\max}\geq1$, we have $\frac{\alpha}{16\sigma}\epsilon^2\geq\frac{\alpha\epsilon^2}{16\sigma_{\max}^2}$ and $\alpha\|s\|^2>\frac{\alpha\epsilon^2}{16\sigma_{\max}^2}$. Hence, by \eqref{eq:teste-misto},
$$
f_i(x^{k+1})\leq f_i(x^k)-\frac{\alpha\epsilon^2}{16\sigma_{\max}^2}=f_i(x^k)-\frac{\alpha}{16c^2}\,\epsilon^{2-2(\beta-1)/\beta}=f_i(x^k)-\frac{\alpha}{\hat\eta}\,\epsilon^{2/\beta},
$$
and the bound follows as in the proof of Theorem~\ref{teo:iteracoes}.

\item[(b)] By the assumptions of item~(b), $\sigma_\ell^k\geq(B+\theta)/2$ for all $\ell$ and $k$. Since $\sigma$ only increases within an iteration, every $\sigma$ used satisfies $\sigma\geq(B+\theta)/2$. Thus, Lemma~\ref{lem:cota-passo}(b) gives $\sigma\|s^k_{(j)}\|\geq G_j(x^k)/4>\epsilon/4$, so that
$$
\alpha\sigma\|s^k_{(j)}\|^2=\frac{\alpha(\sigma\|s^k_{(j)}\|)^2}{\sigma}\geq\frac{\alpha}{16\sigma}\,\epsilon^2,
$$
and the proofs of Lemma~\ref{lem:decr-uniforme}(b) and Theorem~\ref{teo:iteracoes} apply without change.
\end{enumerate}
The proofs of items~(a)--(c) of Theorem~\ref{teo:avaliacoes} only use the number of iterations, Lemma~\ref{lem:decr-uniforme}(a) and the update rule in Step~6, and not the acceptance test. Hence, they hold with $K_\epsilon$ replaced by the corresponding bound.
\end{proof}
\end{appendices}


\begin{thebibliography}{99}

\bibitem{Ama25} Amaral, V.S.: A partially derivative-free cyclic block coordinate descent method for nonseparable composite optimization. Math. Model. Anal. 30, 535--552 (2025). \url{https://doi.org/10.3846/mma.2025.23064}

\bibitem{AABMM22} Amaral, V.S., Andreani, R., Birgin, E.G., Marcondes, D.S., Martínez, J.M.: On complexity and convergence of high-order coordinate descent algorithms for smooth nonconvex box-constrained minimization. J. Glob. Optim. 84, 527--561 (2022). \url{https://doi.org/10.1007/s10898-022-01168-6}

\bibitem{AASS26} Amaral, V.S., Andreani, R., Secchin, L.D., Silva, G.N.: A flexible block coordinate descent method for unconstrained optimization under H\"older continuity. Comput. Optim. Appl. (2026). \url{https://doi.org/10.1007/s10589-026-00827-8}

\bibitem{AAS25} Amaral, V.S., Assunção, P.B., Souza, D.R.: A partially derivative-free proximal method for composite multiobjective optimization in the H\"older setting. Preprint at \url{https://arxiv.org/abs/2508.20071} (2025)

\bibitem{BT13} Beck, A., Tetruashvili, L.: On the convergence of block coordinate descent type methods. SIAM J. Optim. 23, 2037--2060 (2013). \url{https://doi.org/10.1137/120887679}

\bibitem{BM22} Birgin, E.G., Martínez, J.M.: Block coordinate descent for smooth nonconvex constrained minimization. Comput. Optim. Appl. 83, 1--27 (2022). \url{https://doi.org/10.1007/s10589-022-00389-5}

\bibitem{CDM20} Calderón, L., Diniz-Ehrhardt, M.A., Martínez, J.M.: On high-order model regularization for multiobjective optimization. Optim. Methods Softw. 37, 175--191 (2022). \url{https://doi.org/10.1080/10556788.2020.1719408}

\bibitem{CDHS20} Carmon, Y., Duchi, J.C., Hinder, O., Sidford, A.: Lower bounds for finding stationary points I. Math. Program. 184, 71--120 (2020). \url{https://doi.org/10.1007/s10107-019-01406-y}

\bibitem{CGT17} Cartis, C., Gould, N.I.M., Toint, Ph.L.: Worst-case evaluation complexity of regularization methods for smooth unconstrained optimization using H\"older continuous gradients. Optim. Methods Softw. 32, 1273--1298 (2017). \url{https://doi.org/10.1080/10556788.2016.1268136}

\bibitem{chen2023} Chen, J., Tang, L., Yang, X.: A Barzilai--Borwein descent method for multiobjective optimization problems. Eur. J. Oper. Res. 311, 196--209 (2023). \url{https://doi.org/10.1016/j.ejor.2023.04.022}

\bibitem{CTY25} Chen, W., Tang, L., Yang, X.: Generalized conditional gradient methods for multiobjective composite optimization problems with H\"older condition. J. Optim. Theory Appl. 206, 72 (2025). \url{https://doi.org/10.1007/s10957-025-02737-x}

\bibitem{Cla83} Clarke, F.H.: Optimization and Nonsmooth Analysis. Wiley, New York (1983)

\bibitem{das1997} Das, I., Dennis, J.E.: A closer look at drawbacks of minimizing weighted sums of objectives for Pareto set generation in multicriteria optimization problems. Struct. Optim. 14, 63--69 (1997). \url{https://doi.org/10.1007/BF01197559}

\bibitem{Des12} D\'esid\'eri, J.-A.: Multiple-gradient descent algorithm (MGDA) for multiobjective optimization. C. R. Acad. Sci. Paris, Ser. I 350, 313--318 (2012). \url{https://doi.org/10.1016/j.crma.2012.03.014}

\bibitem{FGS09} Fliege, J., Gra\~na Drummond, L.M., Svaiter, B.F.: Newton's method for multiobjective optimization. SIAM J. Optim. 20, 602--626 (2009). \url{https://doi.org/10.1137/08071692X}

\bibitem{FS00} Fliege, J., Svaiter, B.F.: Steepest descent methods for multicriteria optimization. Math. Methods Oper. Res. 51, 479--494 (2000). \url{https://doi.org/10.1007/s001860000043}

\bibitem{FVV19} Fliege, J., Vaz, A.I.F., Vicente, L.N.: Complexity of gradient descent for multiobjective optimization. Optim. Methods Softw. 34, 949--959 (2019). \url{https://doi.org/10.1080/10556788.2018.1510928}

\bibitem{GS05} Gra\~na Drummond, L.M., Svaiter, B.F.: A steepest descent method for vector optimization. J. Comput. Appl. Math. 175, 395--414 (2005). \url{https://doi.org/10.1016/j.cam.2004.06.018}

\bibitem{LP18} Lucambio Pérez, L.R., Prudente, L.F.: Nonlinear conjugate gradient methods for vector optimization. SIAM J. Optim. 28, 2690--2720 (2018). \url{https://doi.org/10.1137/17M1126588}

\bibitem{MG25} Maia, L.F., Gutman, D.H.: The randomized block coordinate descent method in the H\"older smooth setting. Optim. Lett. 19, 1--20 (2025). \url{https://doi.org/10.1007/s11590-024-02161-6}

\bibitem{morovati2016} Morovati, V., Pourkarimi, L., Basirzadeh, H.: Barzilai and Borwein's method for multiobjective optimization problems. Numer. Algorithms 72, 539--604 (2016). \url{https://doi.org/10.1007/s11075-015-0058-7}

\bibitem{Nes12} Nesterov, Yu.: Efficiency of coordinate descent methods on huge-scale optimization problems. SIAM J. Optim. 22, 341--362 (2012). \url{https://doi.org/10.1137/100802001}

\bibitem{Nes15} Nesterov, Yu.: Universal gradient methods for convex optimization problems. Math. Program. 152, 381--404 (2015). \url{https://doi.org/10.1007/s10107-014-0790-0}

\bibitem{PG25} Pinheiro, M.E., Grapiglia, G.N.: Universal nonmonotone line search method for nonconvex multiobjective optimization problems with convex constraints. Comput. Appl. Math. 44, 56 (2025). \url{https://doi.org/10.1007/s40314-024-03018-5}

\bibitem{Pov14} Povalej, Ž.: Quasi-Newton's method for multiobjective optimization. J. Comput. Appl. Math. 255, 765--777 (2014). \url{https://doi.org/10.1016/j.cam.2013.06.045}

\bibitem{Sam26} Sampaio, P.R.: Complexity bounds for smooth multiobjective optimization. J. Complex. 95, 102044 (2026). \url{https://doi.org/10.1016/j.jco.2026.102044}

\bibitem{SK18} Sener, O., Koltun, V.: Multi-task learning as multi-objective optimization. In: Advances in Neural Information Processing Systems 31 (NeurIPS 2018), pp. 527--538 (2018)

\bibitem{SWM25} Soriano, E., Mitkovski, M., Wiecek, M.M.: Decision space decomposition for multiobjective optimization. J. Optim. Theory Appl. 208, 103 (2026). \url{https://doi.org/10.1007/s10957-025-02924-w}

\bibitem{TFY19} Tanabe, H., Fukuda, E.H., Yamashita, N.: Proximal gradient methods for multiobjective optimization and their applications. Comput. Optim. Appl. 72, 339--361 (2019). \url{https://doi.org/10.1007/s10589-018-0043-x}

\bibitem{TV26} Tran, T.H., Vicente, L.N.: Stochastic block coordinate and function alternation for multi-objective optimization and learning. Preprint at \url{https://arxiv.org/abs/2605.12432} (2026)

\bibitem{Tse01} Tseng, P.: Convergence of a block coordinate descent method for nondifferentiable minimization. J. Optim. Theory Appl. 109, 475--494 (2001). \url{https://doi.org/10.1023/A:1017501703105}

\bibitem{War63} Warga, J.: Minimizing certain convex functions. J. Soc. Ind. Appl. Math. 11, 588--593 (1963). \url{https://doi.org/10.1137/0111043}

\bibitem{Wri15} Wright, S.J.: Coordinate descent algorithms. Math. Program. 151, 3--34 (2015). \url{https://doi.org/10.1007/s10107-015-0892-3}

\end{thebibliography}

\end{document}